\documentstyle[amsfonts,amssymb,12pt]{article}
\input xypic
\input cyracc.def

\newfam\cyifam
\font\tencyi=wncyi10
\font\sevencyi=wncyi7
\font\fivecyi=wncyi5

\textfont\cyifam=\tencyi \scriptfont\cyifam=\sevencyi
  \scriptscriptfont\cyifam=\fivecyi
\newfam\cyrfam
\font\tencyr=wncyr10
\font\sevencyr=wncyr7
\font\fivecyr=wncyr5
\def\cyr{\fam\cyrfam\tencyr\cyracc}
\textfont\cyrfam=\tencyr \scriptfont\cyrfam=\sevencyr
  \scriptscriptfont\cyrfam=\fivecyr

\renewcommand{\paragraph}{}

\newcommand{\lon}{\longrightarrow}
\newcommand{\ad}{{\mathrm ad}}
\newcommand{\rar}{\rightarrow}

\newcommand{\Proof}{{\bf Proof}.\, }

\newcommand{\End}{{\mathrm End}}

\newcommand{\Def}{\mathsf{Def}}
\newcommand{\Df}{\mbox{\cyr Def}}

\newcommand{\ovr}{\overrightarrow}
\newcommand{\ovl}{\overline}

\newcommand{\Z}{{\Bbb Z}}
\newcommand{\p}{{\partial}}

\newcommand{\C}{{\Bbb C}}
\newcommand{\R}{{\Bbb R}}
\newcommand{\ot}{\otimes}
\newcommand{\tl}{\tilde}
\newcommand{\Id}{\mbox{Id}}
\newcommand{\Beq}{\begin{equation}}
\newcommand{\Eeq}{\end{equation}}
\newcommand{\Beqr}{\begin{eqnarray}}
\newcommand{\Eeqr}{\end{eqnarray}}
\newcommand{\Beqrn}{\begin{eqnarray*}}
\newcommand{\Eeqrn}{\end{eqnarray*}}
\newcommand{\Ba}{\begin{array}}
\newcommand{\Ea}{\end{array}}
\newcommand{\Bi}{\begin{itemize}}
\newcommand{\Ei}{\end{itemize}}
\newcommand{\Bc}{\begin{center}}
\newcommand{\Ec}{\end{center}}
\newcommand{\fg}{{\frak g}}
\newcommand{\fh}{{\frak h}}
\newcommand{\f}{{\cal O}}

\newcommand{\cM}{{\cal M}}
\newcommand{\cA}{{\cal A}}
\newcommand{\cB}{{\cal B}}

\newcommand{\al}{\alpha}
\newcommand{\be}{\beta}
\newcommand{\ga}{\gamma}
\newcommand{\Ga}{\Gamma}
\newcommand{\Gai}{{\mit{\Upsilon}}}
\newcommand{\var}{\varepsilon}

\newcommand{\tmu}{\tilde{\mu}}

\newcommand{\tv}{\tilde{v}}

\newcommand{\tll}{\tilde{l}}
\newcommand{\tlo}{\tilde{0}}
\newcommand{\tln}{\tilde{1}}
\newcommand{\tk}{\tilde{k}}
\newcommand{\tX}{\tilde{X}}

\newcommand{\tj}{\tilde{j}}

\newcommand{\Ker}{{\mathsf Ker}\, }
\newcommand{\Img}{{\mathsf Im}\, }

\newcommand{\Hom}{\mbox{Hom}}
\newcommand{\bH}{{\bf H}}

\newcommand{\cT}{{\cal T}}
\newcommand{\T}{{\cal T}_{\cal H}}
\newcommand{\X}{{\mbox{\small{\em X}}}}
\newcommand{\Y}{{\mbox{\small{\em Y}}}}

\newcommand{\M}{{\cal H}}

\newcommand{\Ho}{{\cal H}_{\bH}}
\newcommand{\Dg}{D^{\Gamma}}
\newcommand{\dg}{d^{\Gamma}}
\newcommand{\dv}{\delta}

\newcommand{\sip}{\smallskip}
\newcommand{\bip}{\bigskip}

\newcommand{\vst}{\vspace{2 mm}}

\newcommand{\vse}{\vspace{8 mm}}

\begin{document}

\title{Frobenius$_{\infty}$ invariants of homotopy\\ 
Gerstenhaber algebras, I} 
\author{ S.A.\ Merkulov}
\date{}
\maketitle

\begin{abstract}
We construct a functor 
from the derived category of homotopy Gerstenhaber algebras, $\fg$,
with finite-dimensional cohomology  to the purely geometric category of 
so-called $F_{\infty}$-manifolds. The latter contains Frobenius 
manifolds as a  subcategory (so that a pointed Frobenius manifold is itself
 a homotopy Gerstenhaber algebra).
If $\fg$ happens to be formal as a $L_{\infty}$-algebra, then its
$F_{\infty}$-manifold comes equipped with the 
Gauss-Manin connection.
Mirror Symmetry implications are discussed.
\end{abstract}

\section{Introduction} 
Frobenius manifolds play a central
role in the usual formulation of Mirror Symmetry,
as may be seen in the following diagram,
$$
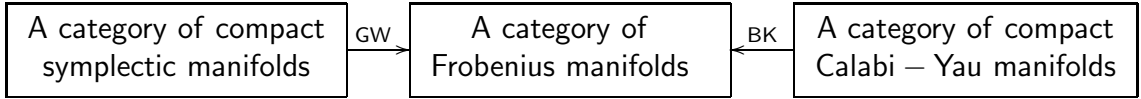
\diagram
{\Ba{c} \mathsf{A\ category\ of\ compact}\\
\mathsf{symplectic\ manifolds}\Ea} \framed \rto^{\mathsf\  GW \ } &
{\Ba{c} \mathsf{A\ category\ of}\\
\mathsf{Frobenius\ manifolds}\Ea} \framed \ \  &
{\Ba{c}\mathsf{A\ category\ of\ compact}\\
\mathsf{ Calabi-Yau\ manifolds}\Ea}\framed  \lto_{\mathsf\ BK\ \ }  
\enddiagram
$$
where morphisms in all categories are just diffeomorphisms preserving
relevant structures, and ${\mathsf GW}$ and ${\mathsf BK}$ stand, 
respectively, for the Gromov-Witten (see, e.g., \cite{Ma})
and Barannikov-Kontsevich \cite{BaKo,Ba} functors. A pair $(\widetilde{M},M)$
consisting of a symplectic manifold  $\widetilde{M}$
 and a Calabi-Yau manifold $M$ is said to be  {\em Mirror}\, if 
${\mathsf GW}(\widetilde{M})={\mathsf BK}(M)$. 
According to Kontsevich \cite{Ko2}, this equivalence is a 
shadow of a more fundamental equivalence of natural $A_{\infty}$-categories
attached to $M$ and $\widetilde{M}$.

\sip

This paper is much motivated by the Barannikov-Kontsevich construction
\cite{BaKo,Ba} of the functor from the right in the above diagram, 
and by Manin's comments \cite{Ma0} on  their construction. The roots 
of the ${\mathsf BK}$ functor lie in the extended 
deformation theory of complex 
structures on  ${M}$, more precisely in very special properties of
the (differential) Gerstenhaber algebra $\fg$ 
``controlling'' such deformations. One of
the miracle features of Calabi-Yau manifolds, the one which played
a key role in the ${\mathsf BK}$ construction, is that deformations 
of their complex 
structures  are  non-obstructed, always producing a {\em smooth}\, versal
 moduli space\footnote{A similar phenomenon occurs in the extended 
deformation 
theory of Lefschetz symplectic structures which also produces, via
the same ${\mathsf BK}$ functor, Frobenius manifolds \cite{Me}.  
These should not
be confused with GW.}. 
In the language of  Gerstenhaber algebras, the exceptional algebraic
properties   necessary to produce a Frobenius manifold out of $\fg$,
have been axiomatized in \cite{Ma0,Ma}. As a result, the functor
$$
\diagram
{\Ba{c} \mbox{``$\mathsf{Exceptional}$''}\\
\mathsf{Gerstenhaber\ algebras}\Ea} \framed & \stackrel{\mathsf BK}{\lon}
&
{\Ba{c} \mathsf{Frobenius}\\
\mathsf{manifolds}\Ea} \framed 
\enddiagram
$$
is now well understood.

\sip

One of our purposes in this paper is to extend the BK functor
from the category of Calabi-Yau manifolds to the category of arbitrary
compact complex manifolds. Which means the study of a diagram
$$
\diagram
{\Ba{c} \mathsf{Arbitrary}\\
\mathsf{Gerstenhaber\ algebras}\Ea} \framed & \stackrel{\mathsf \, ?}{\lon}
&
{\Ba{c} \ \ \ {\mathsf ?} \ \ \ \\
\ \ \  \Ea} \framed 
\enddiagram
$$
Generically, the extended deformation theory of complex structures 
  is obstructed, and it would be  naive to expect that
the question mark above stands for the category of Frobenius manifolds.
In fact, it is not, and the answer is captured in the following notion.

\bip

\paragraph{\bf 1.1. Definition. } An $F_{\infty}$-{\em manifold}\, is
the data $(\M, E,\p, [\mu_{*}], e)$, where
\begin{itemize}
\item[(i)] $\M$ is a formal pointed $\Z$-graded manifold, \\
 \item[(ii)] $E$ is the Euler vector field on $\M$, 
$Ef:=\frac{1}{2}|f|f$, for all homogeneous functions on $\M$
of degree $|f|$,
\item[(iii)] $\p$ is an odd homological (i.e.\ $\p^2=0$) vector field
on $\M$ such that $[E,\p]=\frac{1}{2}\p$ and  
$\p I\subset I^2$, $I$ being the ideal of the 
distinguished point in $\M$, \\
\item[(iv)] $[\mu_n:\otimes^n \T \rar \T]$, $n\in {\Bbb N}$, is a 
homotopy class of smooth unital strong homotopy
commutative ($C_{\infty}$) algebras  defined on the 
tangent sheaf, $\T$, to $\M$, such that
$Lie_E \mu_n= \frac{1}{2} n\mu_n$, for all 
$n\in {\Bbb N}$, and  $\mu_1$ is given by
$$
\Ba{rccc}
\mu_1: & \T & \lon & \T \\
& X &\lon & \mu_1(X):=[\p,X].
\Ea
$$
\\
\item[(v)] $e$ is the unit, i.e.\ 
an even vector field on $\M$ such that $[\p, e]=0$,
$\mu_2(e,X)=X$, $\forall X\in \T$, and $\mu_n(\ldots, e,\ldots )=0$
for all $n\geq 3$.\\
\end{itemize}

\sip

 Clearly, the category of $F_{\infty}$-manifolds containes
Frobenius manifolds as a subcategory.

 \sip

On any $F_{\infty}$-manifold 
the vector field $\p$ defines an integrable distribition, $\Img \mu_1$,
 which is tangent to its subspace
of zeros, $``{\mathsf zeros}(\p)\mbox{''}$.
The structures (i)-(v)  make  the tangent sheaf to
the  smooth part
of the associated quotient, $``{\mathsf zeros}({\p})\mbox{''}/ 
\Img\mu_1$,  into  a sheaf of 
graded unital associative algebras.

\bip

\paragraph{\bf Theorem A. } {\em For any differential 
unital (graded commutative) 
Gerstenhaber
algebra $\fg$, its cohomology, $\bH(\fg)$, if finite-dimensinal, 
is canonically an $F_{\infty}$-manifold.}

\bip

The resulting diagram
$$
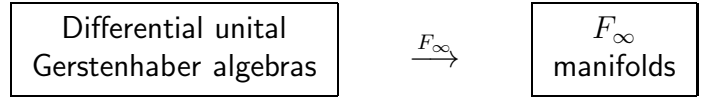
\diagram
{\Ba{c} \mathsf{Differential\ unital}\\
\mathsf{Gerstenhaber\ algebras}\Ea} \framed & \stackrel{F_{\infty}}
{\lon}
&
{\Ba{c} F_{\infty}\\  {\mathsf manifolds }\Ea} \framed 
\enddiagram
$$
implies, in turn, a diagram,
$$
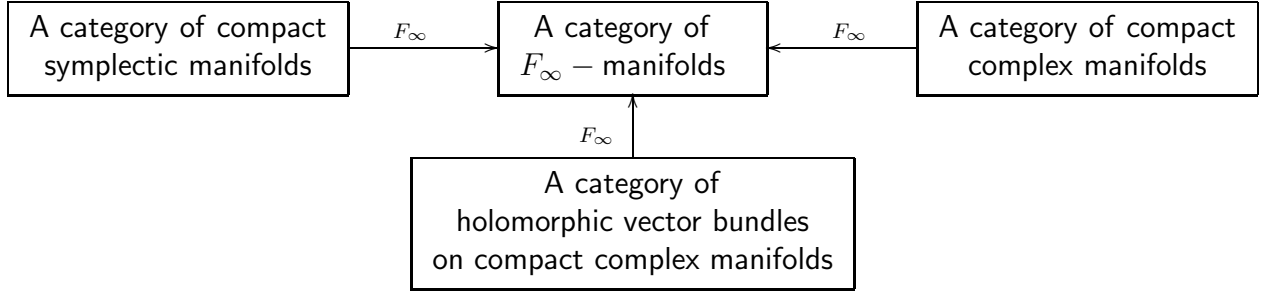
\diagram
{\Ba{c} \mathsf{A\ category\ of\ compact}\\
\mathsf{symplectic\ manifolds}\Ea} \framed \rto^{\  F_{\infty}  } &
{\Ba{c} \mathsf{A\ category\ of}\\
F_{\infty}-\mathsf{manifolds}\Ea} \framed \ \  &
{\Ba{c}\mathsf{A\ category\ of\ compact}\\
\mathsf{complex\ manifolds}\Ea}\framed  \lto_{ F_{\infty}\ \ }  \\
& {\Ba{c} \mathsf{A\ category\ of}\\
\mathsf{holomorphic\ vector\ bundles}\\ \mathsf{on\ compact\ 
complex\ manifolds}\Ea} \framed \uto^{  F_{\infty} \ } &
\enddiagram
$$
through the Gerstenhaber algebras
controlling extended deformations of symplectic, complex and 
holomorphic vector bundle structures. Moreover, the $F_{\infty}$-functor
enjoys the correct ``classical limit'': when restricted
to exceptional Gerstenhaber algebras (i.e.\ the ones
satisfying Manin's axioms \cite{Ma0}), the $F_{\infty}$-functor
coincides precisely with the
Barannikov-Kontsevich  
construction \cite{BaKo}, and hence takes values in the subcategory of 
Frobenius manifolds. 

\sip

Let us emphasize once again that, in the above diagram, 
$F_{\infty}(\mathsf symplectic\ manifolds)$ has nothing to do 
with  Gromov-Witten
invariants\footnote{At best, this is  a very weak shadow of the Mirror
Symmetry, see below.}.
Nevertheless, to rather different mathematical objects
we can canonically attach invariants lying in one and the same geometric
category. Hence we can use these $F_{\infty}$-invariants for a
classification, and even speak about {\em dull mirror symmetry}\, when
$$
F_{\infty}({\mathsf Object}) = F_{\infty}(\widetilde{\mathsf Object}).
$$
Such a relation may be  a shadow of something conceptually
more interesting (cf.\ \cite{Ko2}). 

\sip

Theorem ~A is explained and  generalized by the following

\bip

\paragraph{$\mathbf Theorem\ B$. } {\em 
There is a canonical functor, $F_{\infty}$, from
the  derived category of unital homotopy Gerstenhaber ($G_{\infty}$) 
algebras with finite-dimensional cohomology to the category of
$F_{\infty}$-manifolds.}

\sip

This result implies that the cohomology space of any homotopy Gerstenhaber
algebra is, if finite-dimensional,  canonically an $F_{\infty}$-manifold,
$$
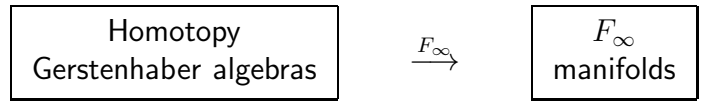
\diagram
{\Ba{c} \mathsf{Homotopy}\\
\mathsf{Gerstenhaber\ algebras}\Ea} \framed & \stackrel{F_{\infty}}
{\lon}
&
{\Ba{c} F_{\infty}\\  {\mathsf manifolds }\Ea} \framed 
\enddiagram
$$

\sip

The recent proof of Deligne's conjecture \cite{Ta1,Ko3,V,MS}
gives the following diagrammatic corollary of Theorem~B,
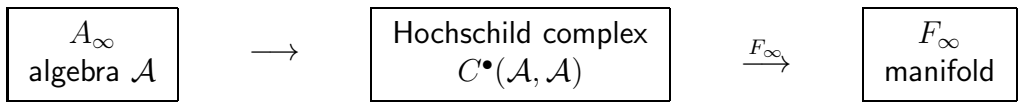
$$
\diagram
{\Ba{c} A_{\infty}\\ {\mathsf algebra}\ {\cA}\Ea} \framed & \lon &
{\Ba{c} \mathsf{Hochschild\ complex}\\
C^{\bullet}(\cA,\cA)\Ea} \framed & 
\stackrel{F_{\infty}}{\lon}
&
{\Ba{c} F_{\infty} \\ {\mathsf manifold} \Ea} \framed 
\enddiagram
$$
which, probably, has  a direct relevance to the Mirror symmetry
through the following specializations, C and D, of statement~B.

\bip

Any $G_{\infty}$-algebra $\fg$ is, in particular, a $L_{\infty}$-algebra
so that its cohomology, $\bH(\fg)$, has the induced structure, 
$[\ , \ ]_{\mathrm ind}$, of Lie algebra. If there exists
a quasi-isomorphism of $L_{\infty}$-algebras, 
$$
\left(\fg, L_{\infty}\hspace{-2mm}-\hspace{-1mm}
{\mathrm  component\ of\ the}\
G_{\infty}\hspace{-2mm}-\hspace{-1mm}{\mathrm structure}
\right) \stackrel{F}{\lon}
\left(\bH(\fg), [\ ,\ ]_{\mathrm ind}\right),
$$
then $\fg$ is said to be $L_{\infty}$-{\em formal},
and $F$ is called a {\em formality map}.
In terms of the associated $F_{\infty}$-invariant, 
the $L_{\infty}$-formality of a $G_{\infty}$-algebra $\fg$ 
gets translated into  a 
canonical flat structure, the Gauss-Manin connection
$\nabla$, 
such that $\nabla_X e=0$ and  
$\nabla_X\nabla_Y\nabla_Z \p=0$ for any flat
vector fields $X,Y$ and $Z$. This specialization of $F_{\infty}$-structure
is called {\em pre}-${\mathit Frobenius}_{\infty}$ structure. 

\bip

\paragraph{\bf Theorem C. } {\em There is a canonical functor
from the category of pairs $(\fg, F)$, where $\fg$ is a
$L_{\infty}$-formal homotopy Gerstenhaber algebra and $F$ a formality
map, to the category of  pre-${\mathit Frobenius}_{\infty}$ 
manifolds.}

\bip

In fact, a pre-${\mathrm Frobenius}_{\infty}$ manifold 
(a Frobenius manifold,
in particular) is itself a homotopy Gerstenhaber algebra. 

\sip

According to Kontsevich's celebrated Formality Theorem \cite{Ko1}, 
for any compact complex manifold
$M$, the Hochschild differential Lie algebra, $C^{\bullet}(\cA,\cA)$,
associated to the  algebra of Dolbeault forms, 
$\cA=(\Gamma(M,\Omega^{0,\bullet}_M), \bar{\p})$,
is $L_{\infty}$-formal so that Theorem~C has a wide area of applications.

\bip

\paragraph{\bf Theorem D.} {\em If a homotopy Gerstenhaber algebra
$\fg$ is quasi-isomorphic, as a $L_{\infty}$-algebra, to an Abelian
differential Lie algebra, then the tangent sheaf, $\T$, to its 
cohomology $\bH(\fg)$
viewed as a linear supermanifold is canonically a sheaf
of unital graded commutative associative algebras.}

\bip

The point is that the Hochschild complex built out of 
the Dolbeault algebra, $\cA=(\Gamma(\cM,\Omega^{0,\bullet}_M), \bar{\p})$,
of a Calabi-Yau manifold $\cM$, satisfies the conditions of Theorem~D.
In fact, the canonically induced associative product on the tangent
sheaf to the associated linear supermanifold, $\bH^{\bullet}(\cM,
\wedge^{\bullet} T_{\cM})$, is, for an appropriate formality map,  
potential and satisfies the WDVV equations. The resulting composition,
$$
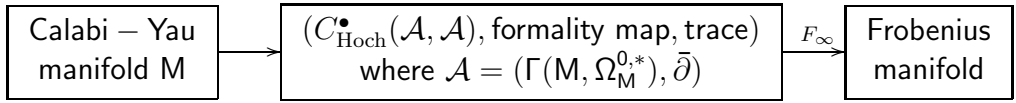
\diagram
{\Ba{c} \mathsf{Calabi-Yau}\\ {\mathsf manifold\ {M}}\Ea} \framed \rto &
{\Ba{c} \left( C^{\bullet}_{\mathrm Hoch}(\cA,\cA), {\mathsf formality\ map},
{\mathsf trace}\right)\\
\mathsf{where}\ \cA=(\Gamma(M,\Omega^{0,*}_{M}), \bar{\p})\Ea} \framed
\rto^{\ \ \ \ \ \ \ \ \ \ \ \ \ \ \ \ \   F_{\infty}} & 
{\Ba{c} {\mathsf Frobenius} \\ {\mathsf manifold} \Ea} \framed 
\enddiagram
$$
together with its analogue for the 
 de Rham algebra of a compact Lefschetz symplectic manifold, 
will be discussed in the second part of this paper.

\sip

\bip

The paper is organized as follows:
\Bi
\item[Section 2:] the origin of the data (i)-(iii) in 
Definition~1.1 of an $F_{\infty}$-manifold is explained via the 
deformation theory. Here we use only
the $L_{\infty}$-component of a $G_{\infty}$-structure. The main 
technical tool is a modified version of the classical deformation functor
which is proved to be non-obstructed. This part of the story is, probably, 
of independent interest.

\item[Section 3:] we use a homotopy techique to explain the origin
of the data (iv)-(v) in Definition~1.1, and to prove Theorems A-D.

\item[Section 4:] we give second proofs of the main claims of this paper
using perturbative solutions of  algebro-differential equations.

\Ei

\section{Deformation functors}

\ \ \ \ \ \  \paragraph{\bf 2.1. Odd Lie superalgebras. } 
Let $k$ be a field with characteristic $\neq 2$. 
An {\em odd Lie superalgebra over $k$}\, is a vector superspace $\fg=\fg_{\tlo}\oplus
\fg_{\tln}$ equipped with an odd $k$-linear map \cite{Ma}
$$
\begin{array}{rccc}
 [\, \bullet\,]: & \fg\ot \fg & \lon & \fg \\
& a \ot b & \lon & [a\bullet b],
\end{array}
$$
which satisfies the following conditions
\Bi
\item[(a)] odd skew-symmetry: $[a\bullet b]=-(-1)^{(\tl{a}+1)(\tl{b}+1)}[b \bullet a]$,
\item[(b)] odd Jacobi identity: 
$$
[a\bullet [b\bullet c]]= [[a\bullet b]\bullet c] + 
(-1)^{(\tl{a}+1)(\tl{b}+1)} [b\bullet [a\bullet c]],
$$
\Ei
for all $a,b,c\in \fg_{\tlo}\cup \fg_{\tln}$.

\sip

The parity change functor transforms this structure into the usual Lie
superalgebra brackets,
$[\ ,\ ]$, on $\Pi\fg$. Thus odd Lie superlgebras are nothing but 
Lie superalgebras in the ``awkwardly'' chosen $\Z_2$-grading.  For 
this reason we  sometimes omit the prefix
{\em odd}, and treat $(\fg,d, [\ \bullet\ ])$ and
$(\Pi\fg, d, [\ ,\ ])$ as different representations of one and the same object.
One advantage of of working with $[\, \bullet \, ]$
rather than with the usual Lie brackets $[\ ,\ ]$ 
is that this awkward $\Z_2$ grading induces the correct, for our purposes,
$\Z_2$-grading on the associated cohomology supermanifold (see below).
Another advantage will become clear below, 
when we introduce
on $\fg$ one more algebraic structure (an {\em even}, in 
this awkward $\Z_2$-grading, 
associative product)  
making $\fg$ into a Gerstenhaber algebra (cf.\ \cite{Ma0,Ma}).

\bip

\paragraph{\bf 2.2. Cohomology as a formal supermanifold.}
A data $(\fg,[\ \bullet \ ],d)$ 
with $(\fg, [\ \bullet \ ])$                          
being a Lie superalgebra and
$$
d:\fg\lon \fg
$$
an odd $k$-linear map satisfying 
\Beqrn
d[a\bullet b] &=& [da\bullet b] - (-1)^{\tl{a}}[a\bullet db] \\
\Eeqrn
is called a {\em differential}\, Lie superalgebra, or shortly
 {\em dLie-algebra}.
 This triple 
$(\fg, [\, \bullet \, ],
 d)$
is often abbreviated to $\fg$.

\sip

The cohomology of  $\fg$,
$$
\bH(\fg):= \Ker d/\Img d\, ,
$$
inherits the structure of Lie superalgebra. We always assume in this paper
that $\bH(\fg)$, which we often abbreviate to $\bH$, 
is a finite dimensional superspace, say $\dim \bH=p|q$. 
Let $\{[e_i], i=1,\ldots,p+q\}$ be a basis
consisting of homogeneous elements with parity denoted by $\tl{i}$,
and $\{t^i, i=1,\ldots,p+q\}$ the associated dual basis in $\bH^*$. 
The supercommutative ring, $k[[t^1,\ldots,t^{p+q}]]$,  of formal 
power series will be abbreviated to $k[[t]]$. 
The (purely) notational advantage of working 
with $k[[t]]$ rather than
with the invariantly defined object $\odot^{\bullet} \bH^*$ 
is that we shall want viewing
\Bi
\item[(i)] $\bH$ as a smooth formal pointed 
$(p|q)$-dimensional supermanifold 
denoted (to emphasize this change of thought) by $\M$ or $\M_{\fg}$,\\
\item[(ii)]  $\{t^i\}$ as linear coordinates on  $\M$, \\
\item[(iii)] $k[[t]]$ as the space of global sections of 
the structure sheaf, $\f_{\M}$,
 on $\M$.
\Ei
The ideal sheaf of the origin, $0\in \M$, will be denote by $I$.

\sip

There is a canonical map
$$
\Ba{rccc}
s:& k[[t]]\ot \bH & \lon & H^0(\M, \T) \vspace{3mm} \\
& \sum a^i(t) [e_i] & \lon & \sum a^i(t) \frac{\p}{\p t^i}, 
\Ea
$$
where $H^0(\M,\T)$ stands for the space of global 
sections of the sheaf, $\T$,
of formal vector fields on $\M$. There is a well defined action
of $H^0(\M,\T)$ 
on both $k[[t]]\ot \bH$ and $k[[t]]\ot \fg$ through the first factor. 
If $\X\,$ is a formal
vector field on $\M$ and $\Gamma$ is an element of 
$k[[t]]\ot \fg$ (or of
$k[[t]]\ot \bH$), then the result of this action is denoted by 
$\ovr{\X\,\,}\,
\Gamma$.

\sip

Any element $\Gamma$ in  $k[[t]]\ot \fg$ (or in  $k[[t]]\ot \bH$)  
can be uniquely decomposed,
$$
\Gamma=\Gamma_{[0]} + \Gamma_{[1]} + \ldots + \Gamma_{[n]} +  \ldots
$$
into homogeneous polynomials, $\Gamma_{[n]}$, of degree $n$ in the variables 
$t^i$. The sum of the first $n$ terms in the above decomposition is 
denoted by $\Gamma_{(n)}$, i.e.\
$\Gamma_{(n)}= \Gamma \, \bmod I^{n+1}$.

\sip

We shall call an element $\Gamma \in k[[t]]\ot \fg$ {\em versal}\,
if it is  even,  $\Gamma \bmod I=0$, 
$\Gamma_{[1]}=\Gamma\bmod I^2\in I\ot \Ker\, d$ 
and
$\Gamma_{[1]}\bmod \Img\, d=\sum_{i=1}^{p+q} t^i [e_i]$.

\sip

The sheaf $\T$ comes canonically equipped with a flat torsion-free
affine connection $\nabla$ whose horizontal sections are, by definition, 
the linear span of $s([e_i])$, $i=1,\ldots, p+q$, i.e.
$$
\nabla \X=0
$$
if and only of $s^{-1}(\X\,)$ is a ``constant'' (independent of $t^i$) 
element in $k[[t]]\ot \bH$. This connection memorizes 
the origin of $\M$ as a vector superspace.

\sip

It will be important, in this paper, to ignore sometimes the 
flat structure and view $\M$ only as a smooth formal supermanifold with a 
distinguished point $0$ but no preferred coordinate system.
To avoid  possible confusion, we adopt from now on this latter 
viewpoint unless the flat connection $\nabla$ is explicitly mentioned.

\bip

\paragraph{\bf 2.2.1. $\Z$-grading.} Differential
 Lie superalgebras $(\fg,d, [\ \bullet\ ])$ 
which we often encounter in geometry
have their $\Z_{2}$-grading induced from a finer structure, $\Z$-{\em grading}, 
which is, by definition, a decomposition of $\fg$ into a direct sum,
$$
\fg=\bigoplus_{i\in \Z} \fg^i,
$$
with the following consistency conditions
\begin{itemize}
\item[(a)]\ \  $d\fg^i \subset \fg^{i+1}$, and\\
\item[(b)]\ \  $[\fg^i\bullet \fg^j]\subset \fg^{i+j-1}$.
\end{itemize}
The $\Z_2$-grading associated to this structure is then simply $\fg_{\tlo}:=\bigoplus_{i\in 2\Z}\fg^i$
and $\fg_{\tln}:=\bigoplus_{i\in 2\Z+1}\fg^i$.

\sip

Clearly, there is an induced $\Z$-grading on the cohomology Lie superalgebra,
$\bH= \oplus_{i\in \Z} \bH^i(\fg)$,  as well as on the structure sheaf, 
$\f_{\M}$, of the associated cohomology supermanifold.

\bip

\paragraph{\bf 2.3. Classical deformation functor.}
One of the approaches to constructing a (versal) deformation space 
of a given mathematical
structure 
$\cA$\, consists of the following  steps (see, e.g. \cite{GM,Ko1,Ba}, 
and references  therein):
\Bi
\item[1)] Associate to $\cA$ a ``controlling" differential $\Z$-graded
Lie algebra $(\fg=\bigoplus_{k\in \Z} \fg^k, d, [\, \bullet\, ])$ over a 
field $k$
(which is usually $\R$ or $\C$).\\

\item[2)] Define the deformation functor
$$
\Ba{rccc}
{\Def}^{\, 0}_{\fg}: & \left\{\Ba{l} \mbox{the category of Artin}\\
                    \mbox{$k$-local algebras} \Ea \right\}&
\lon & \left\{\mbox{the category of sets}\right\}
\Ea
$$
as follows
$$
{\Def}^{0}_{\fg}(\cB)=\left\{\Gamma\in (\fg\ot m_{\cB})^2 \mid d\Gamma
+ \frac{1}{2}[\Gamma\bullet \Gamma]=0\right\}/\exp{(\fg\ot m_{\cB})^1},
$$
where  $m_{\cB}$ is the maximal ideal of the Artin algebra  $\cB$,
the latter is viewed as a $\Z$-graded algebra concentrated in degree 
zero (so that $(\fg\ot m_{\cB})^i=\fg^i\ot m_{\cB}$), and the quotient 
is taken with respect to the following representation of the gauge group 
$\exp{(\fg\ot m_{\cB})^1}$, 
$$
\Gamma \rar \Gamma^g = e^{\ad_g}\Gamma - \frac{e^{ad_g}-1}{\ad_g}dg,
 \ \ \ \ g\in (\fg\ot m_{\cB})^1,
$$
where $\ad$ is just the usual internal automorphism of $\fg$,
$\ad_g \Gamma:= [g\bullet \Gamma]$. \\
\item[3)] Try to represent the deformation 
functor by a topological (pro-Artin) algebra $\f_{\cM}$ so that
$$
{\Def}^{\, 0}_{\fg}(\cB)= \Hom_{\mathrm cont}(\f_{\cM}, \cB).
$$
This associates to the mathematical structure $\cA$ the formal 
 moduli space $\cM$
whose ``ring of functions'' is $\f_{\cM}$.
\Ei

In geometry, one often continues with a fourth step by constructing a 
cohomological splitting of  $\fg$ and applying the Kuranishi
method \cite{Ku2,GM} to represent versally the deformation functor by 
the ring of analytic (rather than formal) functions
on the Kuranishi space.

\sip

The tangent space, ${\Def}_{\fg}^{\, 0}(k[\var]/\var^2)$, to the functor
${\Def}_{\fg}^{\, 0}$ is isomorphic to the  cohomology group 
$\bH^2(\fg)$ of the complex $(\fg,d)$. If one extends in the obvious 
way the above deformation functor to the category of arbitrary 
$\Z$-graded $k$-local Artin algebras (which may not be concentrated 
in degree 0), one gets the functor ${\Def}_{\fg}^{\, *}$
with the tangent space isomorphic to  the full cohomology group
$\oplus _{i\in \Z}\bH^i(\fg)$.

\sip

When working with the extended deformation functor ${\Def}_{\fg}^{*}$
 it is often no loss of essential 
information to forget the $\Z$-grading on $\fg$
and keep only the associated $\Z_2$-grading. One gets then the following
equivalent definition of ${\Def}_{\fg}^{\, *}$:
$$
\Ba{rccc}
{\Def}^{\, *}_{\fg}: & \left\{\Ba{l} \mbox{the category of Artin}\\
                    \mbox{$k$-local superalgebras}\Ea \right\}&
\lon & \left\{\mbox{the category of sets}\right\}
\Ea
$$
$$
{\Def}^{\, *}_{\fg}(\cB):=\left\{\Gamma\in (\fg\ot m_{\cB})_{\tlo} 
\mid d\Gamma
+ \frac{1}{2}[\Gamma\bullet \Gamma]=0\right\}/\exp{(\fg\ot m_{\cB})_{\tln}}.
$$

\sip

This functor is  representable  by 
 a {\em smooth}\, formal moduli space $\cM$
 if there exists a versal (in the sense of Sect.\ 2.2)
solution,
\begin{equation}
\Gamma= \sum_a e_a t^a + \sum_{a_1,a_2}\Gamma_{a_1a_2}t^{a_1}t^{a_2} +
\ldots \  	\in (\fg\ot k[[t]])_{\tlo}, \label{versal}
\end{equation}
to the so-called 
{\em Maurer-Cartan equation},
$$
d\Gamma + \frac{1}{2}[\Gamma\bullet \Gamma]=0.
$$
Due to versality of $\Gamma$, any other
solution over an arbitrary Artin algebra $\cB$ is equivalent to this one
by a  base change $k[[t]]\rar \cB$.

\bip

\paragraph{\bf 2.3.1. $L_{\infty}$-morphisms, part I.}
Let $\fg_1$ and $\fg_2$ be two dLie-algebras.
To formulate the basic theorem of the classical 
deformation theory, we shall need
the following notion: a sequence of linear maps
$$
F_n:\odot^n \fg_1 \lon  \fg_2, \ \ n=1,2, \ldots, \ \ \ \  
\tl{F}_n =\tl{n}+1,
$$
defines a $L_{\infty}$-{\em morphism}\, from $\fg_1$ to $\fg_2$
if 
$$
dF_n(\ga_1, \ga_2,\ldots, \ga_n)
+ \sum_{i=1}^n \pm F_n(\ga_1, \ldots , d\ga_i,\ldots
, \ga_n)
= \ \ \ \ \ \ \ \ \ \ \ \ \ \ \ \ \ \ \ \ \ \ \ \ \ \ \ \ \ \ \ 
\ \ \ \ \ \ 
$$
$$
\ \ \ \ \ \ \ \ \ \ \ \ \ \ \ \ \ \ \ \ \ \
 =\frac{1}{2} \sum_{k+l=n\atop k,l\geq 1}\frac{1}{k!l!}
\sum_{\sigma\in \Sigma_n} \pm \left[
F_k(\ga_{\sigma_1},\ldots , \ga_{\sigma_k}) \bullet
F_l(\ga_{\sigma_{k+1}},\ldots , \ga_{\sigma_n})\right]
+ 
$$
$$
\hspace{-1cm}
+\sum_{i<j}\pm F_{n-1}([\ga_i\bullet \ga_j], \ga_1 , \ldots
, \ga_n),
$$
for arbitrary $\ga_1, \ldots, \ga_n\in \fg_1$.

\sip

In particular, the
first map $F_1$ is a morphism of complexes which 
respects the Lie brackets  up to 
homotopy defined by the second map $F_2$.

\sip

A $L_{\infty}$-morphism 
${F}=\{{F}_n\}:\fg_1\rar\fg_2$
is called a {\em quasi-isomorphism}\, 
if its linear part 
${F}_1$ induces an isomorphism,  $\bH(\fg_1)\rar \bH(\fg_2)$,  of 
associated cohomology groups.

\bip

\paragraph{\bf 2.3.2. Basic Theorem of 
${\Def}$ormation Theory \cite{Ko1} .} {\em 
An $L_{\infty}$-morphism ${F}=\{{F}_n\}:\fg_1\rar\fg_2$
defines a natural transformation of the functors
$$
\Ba{rccc}
{F}_*:& {\Def}^{\, *}_{\fg_1} &\lon& {\Def}^{\, *}_{\fg_2}\vspace{3mm} \\
& \Gamma &\lon & {F}_*(\Gamma):=\sum_{n=1}^{\infty} 
\frac{1}{n!}{F}_n(\Gamma,\ldots,\Gamma),
\Ea
$$
i.e.\ if $\Gamma$ is a  solution to Maurer-Cartan equations  in 
$(\fg_1\ot m_{\cB})_{\tlo}$, then $F_*(\Gamma)$ 
is a solution to Maurer-Cartan 
equations in  $(\fg_2\ot m_{\cB})_{\tlo}$. 

\sip

Moreover, if $F$ is a quasi-isomorphism, then $F_*$
is an isomorphism.
}

\bip

\paragraph{\bf 2.3.3. Corollary.} {\em If a dLie-algebra  
$\fg$ is quasi-isomorphic 
to an Abelian dLie-algebra, then
${\Def}_{\fg}^{\, *}$ is versally representable
by a smooth formal pointed supermanifold $\M_{\fg}$
(the cohomology supermanifold of $\fg$).
}

\bip

\Proof  If $\fh$ is an Abelian dLie-algebra,
then, in the notations of Sect.\ 2.2, 
$$
\Gamma= \sum_{i=1}^{\dim \bH(\fh)}t^i [e_i]
$$ 
is a versal solution of Maurer-Cartan equations.
Hence ${\Def}_{\fh}^{\, *}$ is representable by $\M_{\fh}$.

\sip

If $\fg$ is quasi-isomorphic to $\fh$, 
the required statement follows from Theorem~2.3.2
and the isomorphism $\M_{\fg}=\M_{\fh}$. 
\hfill $\Box$

\bip

There are two remarkable examples, one dealing with extended deformations of 
complex structures on a Calabi-Yau manifold \cite{BaKo} and another
with extended deformations of the symplectic structure on a Lefschetz manifold
\cite{Me}, when the rather strong condition of Corollary~2.3.3
 holds true\footnote{One more example  (of a 
different technical origin though of  the same mirror symmetry flavor)
when a naturally extended deformation problem  gives rise to a {\em smooth}\, 
extended moduli space
is discussed in \cite{Me2}.}. In general, however,
there will be obstructions to constructing a versal solution to the 
Maurer-Cartan equations, and we will have to resort to other technical means
such as modifying the deformation functor
as explained below in Sect.\ 2.4 below or further 
extending the 
deformation problem to the category of $L_{\infty}$-algebras. 

\bip
\paragraph{\bf 2.3.4. Example (deformations of complex manifolds).}
It is well known that the total space of the cotangent bundle, 
$\Omega^1_{\R}$, to a  real
$n$-dimensional manifold 
$M$ carries a natural Poisson structure $\{\, , \, \}$ making the structure 
sheaf
$\f_{\Omega^1_{\R}}$  into a sheaf of Lie algebras. In a natural
local coordinate system $(x^a, p_a:=\p/\p x^a)$, 
$$
\{f,g\}= \frac{\p f}{\p p_a}\frac{\p g}{\p x^a} - \frac{\p f}{\p x^a}
\frac{\p g}{\p p_a}.
$$

\sip 

If we now change the parity of the fibers
of the natural projection $\Omega^1_{\R}\rar M$ (which is allowed 
since they are vector spaces), we will get an $(n|n)$-dimensional
supermanifold $\Pi \Omega^1_{\R}$ equipped with a natural odd Poisson
structure $\{\ \bullet \ \}$ making the structure sheaf
$\f_{\Pi\Omega^1_{\R}}$  into a sheaf of odd Lie superalgebras. In a natural
local coordinate system $(x^a, \psi_a:=\Pi \p/\p x^a)$ on $\Pi \Omega^1_{\R}$, 
$$
\{f\bullet g\}= \frac{\p f}{\p \psi_a}\frac{\p g}{\p x^a} + (-1)^{\tl{g}(
\tl{f}+1)} 
\frac{\p g}{\p \psi_a}
\frac{\p f}{\p x^a}.
$$
The smooth functions on $\Pi\Omega^1_{\R}$ have a simple geometric 
interpretation in term of the underlying manifold $M$ --- they are just 
smooth polyvector fields. Indeed, a standard power series decomposition 
in odd variables gives
$$
f = \sum_{k=0}^n \sum_{a_1,\ldots,a_k} f^{a_1\ldots a_k}(x)\psi_{a_1}\ldots
\psi_{a_k}
$$
implying the isomorphism of sheaves $\f_{\Pi\Omega^1_{\R}}= 
\Lambda^{\bullet}T_{\R}$, where $T_{\R}$ is the real tangent bundle to $M$.

\sip

Therefore, the sheaf $\Lambda^{\bullet} T_{\R}$ with the $\Z_2$-grading,
$$
(\Lambda^* T_{\R})_{\tlo}:=\Lambda^{\mathrm even }T_{\R}, \ \ \ \ \ \ \ \ \ \ \ \
(\Lambda^* T_{\R})_{\tln}:=\Lambda^{\mathrm odd }T_{\R},
$$
 induced from that
on $\f_{\Pi\Omega^1_{\R}}$, 
  is naturally a sheaf of odd Lie superalgebras.
The odd Poisson bracket $\{\, \bullet\ \}$ is called, in this incarnation,
the {\em Schouten}\, bracket and is often denoted by 
$[\, \bullet\ ]_{\mathrm Sch}$. 

\sip

If $M$ is a complex manifold, then the  canonical odd Poison structure on
the parity changed holomorphic cotangent bundle, $\Pi \Omega^1_M$,
is itself holomorphic giving rise thereby to the structure 
of odd Lie superalgebra
on the sheaf, $\Lambda^{\bullet}T_M$, of holomorphic polyvector fields.
This can be used to make the vector space
$$
\fg=\bigoplus_{k\in \Z} \fg^k, \ \ \ \ \ \ \ \ \  
\fg^k :=\bigoplus_{i+j=k}\Gamma(M, 
\Lambda^i T_M \ot \Lambda^j \overline{T}_M^*),
$$
into a $\Z$-graded differential  algebra by taking  $\bar{\p}$, 
the $(0,1)$ part of the de Rham operator, as a differential, and the
map,
$$
\begin{array}{rccc} 
 [\, \bullet\,]: & \Gamma(M,\Lambda^{i_1} T_M \ot 
\Lambda^{j_1}\ovl{T}_M^*)
\, \times \, \Gamma(M,\Lambda^{i_2} T_M \ot \Lambda^{j_2}\ovl{T}_M^*) 
& \rar & \Gamma(M,
\Lambda^{i_1+i_2-1} T_M \ot \Lambda^{j_1+j_2}\ovl{T}_M^*)  \\
& \X_1\ot \ovl{w}_1 \, \times \X_2\ot \ovl{w}_2  & \rar & 
[\X_1\ot \ovl{w}_1\bullet \X_2\ot \ovl{w}_2],
\end{array}
$$
given by
$$
[\X_1\ot \ovl{w}_1\bullet \X_2\ot \ovl{w}_2]:= (-1)^{\tl{j}_1 \tl{i}_2}
[\X_1\bullet \X_2]_{\mathrm Sch}\ot (\ovl{w}_1\wedge \ovl{w}_2).
$$
as the (odd) Lie brackets.

\sip

The importance of this dLie-algebra stems from the fact that
the associated deformation functors ${\Def}^{\, 0}_{\fg}$ and
${\Def}^{\, *}_{\fg}$ describe, respectively, ordinary and
 extended deformations of the given complex structure on a smooth
 manifold $M$. Indeed, a complex structure on a real
$2n$-dimensional manifold $M$ is a decomposition, $\C\ot T_{\R}=
T_M\oplus \ovl{T}_M$, of the complexified real tangent bundle into a 
direct sum of complex integrable distributions, $T_M$ and its complex 
conjugate $\ovl{T}_M$. Another  decomposition like that,
$\C\ot T_{\R}=
T'_M\oplus \ovl{T}'_M$, can be described in terms of the original 
complex structure by the graph, $\ovl{T}_M'$, of a linear map 
$\Gamma: \ovl{T}_M\rar 
{T}_M$, i.e.\ by an element $\Gamma\in \fg^2$. The integrability of $T'_M$
amounts  then to the Maurer-Cartan equation in $\fg^2$,
$$
\ovl{\p} \Gamma + \frac{1}{2}[\Gamma\bullet \Gamma]_{\mathrm Sch}=0.
$$
By solving (if possible) the above equation in the full Lie algebra
$\fg$ rather than in its subalgebra $\fg^2$ (and taking the quotient
by the gauge group describing equivalent deformations),
one gets a so-called extended complex structure on $M$ whose 
geometric meaning is not yet fully understood. It is understood \cite{Ko2}, 
however, 
that this structure does have an important Mirror Symmetry aspect
(at least for Calabi-Yau manifolds):
Kontsevich noticed that his
Formality Theorem \cite{Ko1} identifies the moduli space of 
extended complex structures
on a given complex manifold $M$
with the moduli space of $A_{\infty}$-deformations 
of the derived category of coherent sheaves on $M$, 
mirror counterpart 
 of
the conjectured Fukaya category built out of a dual complex manifold 
$\widetilde{M}$.

\sip

If $M$ is a Calabi-Yau manifold, then, as was shown by Barannikov and 
Kontsevich \cite{BaKo}, 
the  Maurer-Cartan equations in $\fg$ admit
a versal solution of the form (\ref{versal}) implying that the moduli space 
of extended deformations
of complex structures on $M$ is smooth, and is isomorphic\footnote{Strictly 
speaking, this isomorphism holds true in the category of formal manifolds, 
in which we work in this paper. It is no problem to choose the power 
series (\ref{versal}) convergent thereby inducing on
the extended moduli space  a smooth analytic structure. 
The latter is then analytically isomorphic to an open neighbourhood 
of zero in $\bH$ which we denote  by the same symbol 
$\M$ (and continue doing this every time the analyticity aspect emerges).}
to $\M$.  In fact they 
have shown much more \cite{BaKo,Ba}:
in this case $\M$ has an  induced structure of Frobenius manifold
which conjecturally coincides with the Frobenius manifold structure
on $H_*(\widetilde{M},\C)$ constructed via the Gromov-Witten invariants
of the dual Calabi-Yau manifold $\widetilde{M}$.  Barannikov \cite{Ba} 
checked 
this conjecture
for projective complete intersections Calabi-Yau manifolds.

\sip

For a general complex manifold $M$, the (extended) deformations are obstructed
and the Maurer-Cartan equations associated to $\fg$
have no versal solutions of the form (\ref{versal}). 
It is one of the main tasks of this paper
to understand what happens to the Barannikov-Kontsevich's 
Frobenius  structure on $\T$ in 
the presence of obstructions.

\bip

\paragraph{\bf 2.3.5. Example (deformations of Poisson and symplectic 
manifolds).} It is well known that a 2-vector field, 
$\nu_0\in \Gamma(M, \Lambda^2 T_{\R})$, defines a Poisson structure,
$$
\{f,g\} = \nu_0(df\ot dg), \ \ \ \ \ f,g \in \f_M,
$$
on a smooth real $n$-dimensional manifold $M$ if and only if 
$$
[\nu_0\bullet \nu_0]_{\mathrm Sch}=0.
$$
Then a deformed 2-vector field, $\nu_0 + \nu \in \Gamma(M, \Lambda^2 T_{\R})$,
is again a Poisson structure if and only if $\nu$ satisfies
the Maurer-Cartan equation,
$$
d\nu + \frac{1}{2}[\nu\bullet \nu]_{\mathrm Sch}=0,
$$
in the differential $\Z$-graded algebra
$$
\left(\fg=\bigoplus_{i=0}^n \Gamma(M, \Lambda^i T_{\R}),\ 
[\ \bullet\ ]_{\mathrm
Sch}, \
d:=[\nu_0\bullet\ldots ]_{\mathrm Sch}\right).
$$
Hence
the associated deformation functors 
${\Def}^{\, 0}_{\fg}$/${\Def}^{\, *}_{\fg}$ describe
(extended) deformations of the given Poisson structure $\nu_0$
on $M$. 

\sip

For  generic $\nu_0$, the associated cohomology group $\oplus_i\bH^i(\fg)$
may not be finite-dimensional even for compact manifolds. 
If, however, $\nu_0$ comes from a symplectic $2$-form $\omega$ on $M$, 
the situation is very different. In this case one may use 
the ``lowering indices map''  $\omega: T_{\R} \rar \Omega^1_{\R}$ to 
identify $(\fg,d)$ with the de Rham complex
on $M$ --- it is not hard to check that this map sends the differential
$[\nu_0\bullet\ldots]_{\mathrm Sch}$ on $\Lambda^{\bullet} T_{\R}$
into the usual de Rham differential on $\Omega^{\bullet}_{\R}$. 
The (odd) Lie brackets induced on $\Omega^*_{\R}$ from the Schouten 
brackets on $\Lambda^{\bullet} T_{\R}$   we denote
by $[\ \bullet\ ]_{\omega}$ to emphasize its dependence on the 
symplectic structure. Hence the deformation functors
${\Def}^{\, 0}_{\fg}$/${\Def}^{\, *}_{\fg}$ associated with the 
dLie-algebra
$$
\left(\fg=\bigoplus_{i=0}^n \Gamma(M, \Omega^i_{\R}),\  [\ \bullet\ ]_{\omega},
\
d=\mbox{de Rham differential}\right)
$$
describe
(extended) deformations of a symplectic structure $\omega$
on $M$. Its cohomology $\bH$ is nothing but the de Rham cohomology of $M$.

\sip

A compact symplectic manifold $(M,\omega)$ is called {\em Lefschetz} \,
if the  the natural cup product on the de Rham cohomology,
$$
[\omega^k]: H^{m-k}(M,\R) \lon H^{m+k}(M, \R)
$$
is an isomorphism for any $k\leq m:=\frac{1}{2}\dim M$.
This class of manifolds (which includes the class of K\"ahler 
manifolds by the Hard Lefschetz Theorem)
 is of interest to us  for 
the  extended deformation functor ${\Def}^{\, *}_{\fg}$
associated with an arbitrary Lefschetz symplectic manifold is non-obstructed
and 
is  representable by a smooth moduli space isomorphic to 
$\M$;  moreover, this moduli space of ``extended symplectic structures'' 
is always a Frobenius manifold \cite{Me}. This result is more than parallel
to the Barannikov-Kontsevich's construction of extended 
moduli spaces/Frobenius structures
for Calabi-Yau manifolds --- it is just another example when their beautiful
machinery works (see a nice exposition of Manin \cite{Ma}).

\sip

The extended deformation functor associated with a generic
compact symplectic manifold seems to be  obstructed,
and one should employ different techniques (see below) 
to study geometric structures
induced on moduli spaces of extended symplectic structures.

\bip

\paragraph{\bf 2.3.6. Example (deformations of holomorphic vector bundles).} 
Let $E\rar M$ be a holomorphic vector bundle on a complex $n$-dimensional
manifold $M$.
There is an associated differential $\Z$-graded Lie algebra 
$$
\left(\fg=
\Gamma(M, \End E \ot \ovl{\Omega}^{\bullet}_M[-1]), \ 
[\ \bullet\ ],\ \bar{\p} \right) 
$$
with the Lie brackets, 
$$
\begin{array}{rccc} 
 [\, \bullet\,]: & \Gamma(M,\End E \ot \ovl{\Omega}^{i_1-1}_M)
\, \times \, \Gamma(M, \End E\ot \ovl{\Omega}^{i_2-1}_M)
& \rar & 
\Gamma(M,\End E \ot \overline{\Omega}^{(i_1+i_2-1)-1}_M) \\
& A_1\ot \ovl{w}_1 \, \times A_2\ot \ovl{w}_2  & \rar & 
[A_1\ot \ovl{w}_1\bullet A_2\ot \ovl{w}_2],
\end{array}
$$
given by
$$
[A_1\ot \ovl{w}_1\bullet A_2\ot \ovl{w}_2]:= (A_1A_2-A_2A_1)\ot
\left(\ovl{w}_1\wedge \ovl{w}_2\right) .
$$

\bip

The deformation functor ${\Def}^{\, 0}_{\fg}$
associated with this algebra describes
deformations of the holomorphic structure in 
the vector bundle $E$.
It is tempting to view
its extension ${\Def}^{\, *}_{\fg}$ as a tool for studying  {\em extended}\,
deformations, but we reserve this role for the 
functor ${\Def}^{\, *}_{\ldots}$  associated
with a larger differential algebra constructed in 3.1.5 below.

\sip

In general, all these functors are obstructed.

\sip

One may combine this differential Lie algebra (or its extension 3.1.5) 
together with the one of  Example~2.3.4
into their natural semi-direct product
to study {\em joint}\, (extended) deformations of the pair $E\rar M$.

\bip

\paragraph{\bf 2.4. $L_{\infty}$-algebras.}
These algebras will play only an auxiliary, purely technical,
role in this paper.

\sip

By definition, 
a {\em strong homotopy Lie algebra}, or shortly $L_{\infty}$-{\em algebra},
is a vector superspace $\fh$ equipped with linear maps,
$$
\Ba{rcccc}
\mu_k: & \Lambda^k \fh & \lon & \fh & \\
& v_1\wedge\ldots \wedge v_k & \lon & \mu_k(v_1, \ldots, v_k), 
& \ \ \ \ k\geq 1, \ \ \ \  \tl{\mu}_k=\tl{k},
\Ea
$$
satisfying, for any $n\geq 1$ and
 arbitrary $v_1, \ldots, v_n \in \fh_{\tlo}\cup \fh_{\tln}$,
 the following
{\em higher order Jacobi identities},
$$ 
\sum_{k+l=n+1}\sum_{\sigma\in Sh(k,n)} (-1)^{\tilde{\sigma}+k(l-1)} 
e(\sigma;v_1,\ldots,v_n)
 \mu_l\left(\mu_k(v_{\sigma(1)},\ldots,v_{\sigma(k)}),
v_{\sigma(k+l)},\ldots, v_{\sigma(n)}\right)=0,
$$
where $Sh(k,n)$ is the set of all permutations $\sigma:\{1, \ldots, n\} \rar
\{1,\ldots,n\}$ which satisfy $\sigma(1)<\ldots< \sigma(k)$ and 
$\sigma(k+1)<\ldots<\sigma(n)$. The symbol $e(\sigma; v_1, \ldots, v_n)$ 
(which we abbreviate from now on to $e(\sigma)$)
 stands for the {\em Koszul
sign}\,  defined by the equality
$$
v_{\sigma(1)}\wedge \ldots \wedge v_{\sigma(n)}= (-1)^{\tilde{\sigma}}e(\sigma)
v_1\wedge \ldots \wedge v_n,
$$
$\tilde{\sigma}$ being the parity of the permutation
$\sigma$.

\vst

The $\Z$-graded version of this definition would require $\mu_k$
to be homogeneous of degree $2-k$.

\sip

This notion as well as the associated notion of  $A_{\infty}$-algebra 
(reminded below) are due to Stasheff \cite{St1}.

\vst

The first three higher order Jacobi identities have the form
\begin{description}
\item[$n=1$:]\ \  $d^2=0$,\\
\item[$n=2$:]\ \  $d[v_1, v_2]= [dv_1, v_2] + (-1)^{\tv_1} [v_1,
dv_2]$,\\
\item[$n=3$:]\ \ $[[v_1, v_2],v_3] +
(-1)^{(\tv_1+\tv_2)\tv_3}[[v_3, v_1],v_2]
+ (-1)^{\tv_1(\tv_2+\tv_3)}[[v_2, v_3],v_1]$ =
$-d\mu_3(v_1,v_2,v_3) - \mu_3(dv_1,v_2,v_3) - (-1)^{\tv_1}\mu_3(v_1,dv_2,v_3)
-(-1)^{\tv_1+\tv_2}\mu_3(v_1,v_2,dv_3)$,\\
\end{description}
where we denoted $dv_1:=\mu_1(v_1)$ and $[v_1, v_2] := \mu_2(v_1,v_2)$.

\vst

Therefore, $L_{\infty}$-algebras with $\mu_k=0$ for $k\geq 3$
are nothing but the usual  differential Lie superalgebras with the
differential $\mu_1$ and the Lie bracket
$\mu_2$. If, furthermore, $\mu_1=0$, one gets  the class of usual Lie
superalgebras.

\bip

\paragraph{\bf 2.4.1. Odd $L_{\infty}$-algebras.}
To make the above picture consistent with the choices made in  Sect.\ 2.1,
we should change the parity of $\fh$. Hence we shall work
from now on  with $\fg:=\Pi\fh$
and denote 
$$
\mu_n(v_1\bullet v_2\bullet\ldots \bullet v_n):= \Pi \mu_n(\Pi v_1,
\Pi v_2, \ldots, \Pi v_n), \ \ \ \forall v_1, \ldots, v_n \in \fg,
$$
for all $n\geq 1$. This change of grading also unveils, through the 
following three observations, 
a rather compact image of the $L_{\infty}$-structure itself:
\begin{itemize}
\item[1)]  
The vector superspace $\odot^{\bullet}\fg=\bigoplus_{n=1}^{\infty}\odot^n \fg$
has a natural structure of cosymmetric {\em coalgebra},
$$
\Delta(w_1\odot\ldots \odot w_n) =\sum_{i=1}^n \sum_{\sigma\in
Sh(i,n)} e(\sigma) \left(w_{\sigma(1)}\odot \ldots \odot
w_{\sigma(i)}\right)
\otimes \left(w_{\sigma(i+1)}\odot \ldots \odot
w_{\sigma(n)}\right).
$$
\\
\item[2)] Every coderivation of this coalgebra, i.e. an odd map 
$Q:\odot^{\bullet} \fg \rar \odot^{\bullet} \fg$ satisfying 
$\Delta \circ Q=\left(Q\otimes \Id + \Id\otimes Q\right)\circ
\Delta$, is equivalent to an arbitrary series of odd linear maps,
$\mu_n:\odot^n \fg \rar \fg$. \\
\item[3)] A codifferential $Q=\{\mu_*\}$ is a differential, i.e.\
 $Q^2=0$, if and only if $\mu_n$ satisfy the 
higher order Jacobi identities.
\end{itemize}

\sip

In conclusion, an (odd) $L_{\infty}$-structure on $\fg$
is equivalent to a codifferential on $(\odot^{\bullet} \fg, \Delta)$.

\bip

\paragraph{\bf 2.4.2. $L_{\infty}$-morphisms, part II.} Given two
$L_{\infty}$-algebras, $(\fg, \mu_*)$ and $(\fg', \mu'_*)$. A
$L_{\infty}$-morphism $F$ from the first one to the second is, by
definition, a differential coalgebra 
homomorphism
$$
F: \left(\odot^{\bullet}\fg, \Delta, Q\right) \lon \left(\odot^{\bullet} \fg', \Delta,
Q'\right),
$$
i.e.\ a linear map $F: \odot^{\bullet} \fg \rar \odot^{\bullet} \fg'$
satisfying $(F\ot F)\circ \Delta = \delta'\circ F$ and $F\circ Q=Q'\circ F$.
The first of these equations says that $F$
 is completely determined by a set of  linear maps 
${F}_n: \odot^n \fg
\rar \fg'$ of parity $\tl{n} + \tln$, 
while the second one imposes on these $F_n'$ a sequence of 
linear equations. If both input and output of 
$F$ are
usual  differential Lie superalgebras, these equations  are
presicely the ones  written down
in Sect.\ 2.3.1.

\vst

A $L_{\infty}$-morphism $F:(\fg, \mu_*)\rar (\fg', \mu'_*)$ is
called a {\em quasi-isomorphism}\, if its first component 
$F_1: \fg\rar \fg'$
induces an isomorphism between cohomology groups of complexes
$(\fg, \mu_1)$ and $(\fg', \mu'_1)$. 
It is called a {\em
homotopy of the $L_{\infty}$-algebras}, if 
$F_1: \fg\rar \fg'$ is an isomorphism of underlying
vector graded superspaces.

\sip

If the
$L_{\infty}$-morphism $F:(\fg, \mu_*)\rar (\fg', \mu'_*)$ is
 a quasi-isomorphism, then, as was proved in \cite{Ko1}, 
there exists a $L_{\infty}$-morphism
 $F':(\fg', \mu'_*)\rar (\fg, \mu_*)$ which induces the inverse
 isomorphism between cohomology groups of complexes
$(\fg, \mu_1)$ and $(\fg', \mu'_1)$.

\sip

Two $L_{\infty}$-morphisms, $F,G:(\fg, \mu_*)\rar (\fg', \mu'_*)$,
are said to be {\em homotopy equivalent}\, if there is an odd
linear map $h: \odot^{\bullet} \fg \rar \fg'$ such that
$$
\Delta'\circ h = F\ot h + h\ot G, \ \ \ {\mathrm and}\ \ \ \
F=G + Q'\circ h + h\circ Q.
$$
This map is completely
determined by a set of linear maps, $\{h_n: \odot^n \fg \rar \fg',\, 
\tl{h_n}=\tl{n},
n=1,2,\ldots \}$,
and 
is called a {\em homotopy of morphisms}.
The resulting homotopy relation on the set of all $L_{\infty}$-morphisms
from $(\fg, \mu_*)$ to $(\fg', \mu'_*)$ is an equivalence relation.

\bip

\paragraph{\bf 2.4.3. A geometric interpretation of a
$L_{\infty}$-algebra $\fg$.} The dual of the free cocommutative coalgebra
$\odot^{\bullet}\fg$ can be identified with the algebra of
formal power series on the vector superspace $\fg$ viewed as a
formal pointed supermanifold (to emphasize this change of thought 
we denote the supermanifold structure on $\fg$ by $M_{\fg}$). With
this identification, the $L_{\infty}$-structure $\mu_*$ on $\fg$,
that is the  codifferential $Q$ on $\odot^{\bullet} \fg$, goes
into an odd vector field $Q$ on $M_{\fg}$ satisfying \cite{Ko1}
\begin{itemize}
\item[a)] $Q^2=0$,
\item[b)] $Q I\subset I$,
\end{itemize}
where $I$ is the ideal of the distinguished point, $0\in M_{\fg}$.
(An odd vector field on a formal pointed superspace satisfying the above
two conditions is usually called {\em homological}.)

A $L_{\infty}$-morphism $F$ between two $L_{\infty}$-algebras $(\fg,\mu_*)$
and $(\fg',\mu_*')$ is nothing but  a $Q$-equivariant map
between the associated formal pointed homological 
supermanifolds, $(M_{\fg}, Q)$ and
$(M_{\fg'}, Q')$.

\sip

\bip

\paragraph{\bf 2.5. A modified deformation functor.}  For a general
dLie-algebra $\fg$,
the classical deformation functor $\Def_{\fg}^{\, *}$
is not representable by a  {\em smooth}\, versal moduli space. At best
one can use Kuranishi technique to represent  $\Def_{\fg}^{\, *}$ 
by a singular analytic space. There is, however, a simple geometric way
to keep track of  versality and smoothness.
The idea is as follows.

\sip

First, we extend the input category 
used in the construction of ${\Def}_{\fg}^*$ in Sect.\ 2.3
to a  category  of {\em differential}\, Artin superalgebras
whose $\mathsf{Ob}$s are pairs, $(\cB, \p)$, consisting of an Artin 
superalgebra $\cB$ together with a differential 
$\p:\cB\rar\cB$ satisfying $\p m_{\cB}\subset
m_{\cB}^2$, and whose $\mathsf{Mor}$s are morphisms of Artin superalgebras
commuting with the differentials. A $\Z$-graded version of this definition
would require $\p$ to have degree $+1$.

\sip

Second, to the ``controlling'' differential Lie algebra $\fg$
we associate a new deformation functor,  
$$
\Ba{rccc}
{\Df}^{*}_{\fg}: & \left\{\Ba{c} \mbox{the category 
   of differential}\\
                    \mbox{Artin superalgebras}
                    \Ea \right\}&
\lon & \left\{\mbox{the category of sets}\right\}\vspace{3mm}\\
& (\cB,\p) & \lon & {\Df}^{\, *}_{\fg}(\cB,\p)
\Ea
$$
by setting
$$
{\Df}^{ *}_{\fg}(\cB,\p)=
\left\{\Gamma\in (\fg\ot m_{\cB})_{\tlo} \mid d\Gamma
+ \vec{\p}\Gamma + \frac{1}{2}[\Gamma\bullet \Gamma]=0\right\}/
\exp{(\fg\ot m_{\cB})_{\tln}}.
$$
Here the quotient is taken with respect to the 
following representation of the 
gauge group, 
$$
\Gamma \rar \Gamma^g = e^{\ad_g}\Gamma - 
\frac{e^{ad_g}-1}{\ad_g}(d+\vec{\p})g,
 \ \ \ \ \forall g\in (\fg\ot m_{\cB})_{\tln}.
$$

\sip

On the subcategory $(\cB, 0)$ the deformation functor
${\Df}^{*}_{\fg}$ coincides precisely
with the classical one 
$\mathsf{Def}^{\, *}_{\fg}$.

\bip

\paragraph{\bf 2.5.1. Remark.} If, for a derivation $\p:\cB\rar \cB$,
an element
 $\Gamma\in (\fg\ot m_{\cB})_{\tlo}$ satisfies the equation
(which we sometimes call the {\em Master equation}),
$$
 d\Gamma+ \vec{\p}\Gamma + \frac{1}{2}[\Gamma\bullet \Gamma]=0,
$$
then
\Beqrn
0 &=& d\left( d\Gamma+ \vec{\p}\Gamma + \frac{1}{2}[\Gamma\bullet \Gamma]
\right) \\
&=&  -\vec{\p} d\Gamma + [d\Gamma\bullet \Gamma]\\
&=&  -\vec{\p} d\Gamma - [\vec{\p}\Gamma\bullet \Gamma]
-\frac{1}{2}\left[[\Gamma\bullet\Gamma]\bullet 
\Gamma\right] \\
&=& - \vec{\p}\left( d\Gamma + \frac{1}{2}
\left[\Gamma\bullet\Gamma\right]
\right) \\
&=& - \vec{\p}^{\, 2}\Gamma,
\Eeqrn
motivating our assumption above 
that $\p$ is a differential in $\cB$ rather than merely a derivation.

\bip

\paragraph{\bf 2.5.2. $L_{\infty}$-extension of $\Df^*$.}
This  extension   will be used later only as a 
technical tool  in the study of $\Df^*_{\fg}$
for usual  differential Lie superalgebras $\fg$.

\sip

Given a $L_{\infty}$-algebra $(\fg, \mu_*)$,
we define,  
$$
\Ba{rccc}
{\Df}^{*}_{\fg}: & \left\{\Ba{c} \mbox{the category 
   of differential}\\
                    \mbox{Artin superalgebras}
                    \Ea \right\}&
\lon & \left\{\mbox{the category of sets}\right\}\vspace{3mm}\\
& (\cB,\p) & \lon & {\Df}^{\, *}_{\fg}(\cB,\p)
\Ea
$$
by setting
$$
{\Df}^{ *}_{\fg}(\cB,\p)=
\left\{\Gamma\in (\fg\ot m_{\cB})_{\tlo} \mid 
 \vec{\p}\Gamma =
\sum_{n=1}^{\infty} \frac{(-1)^{n(n+1)/2}}{n!}
\mu_n(\Gamma \bullet \ldots \bullet \Gamma)\right\}/\sim
$$
Here the quotient is taken with respect to the 
 gauge equivalence, $\sim$, 
which is best described using
the following geometric model of the $\Df{\mathrm ormation}$ functor.

\sip

In ${\mathsf Category}^{op}$, both the differential Artin superalgebra, 
$(\cB,\p)$, and the $L_{\infty}$-algebra, $(\fg, \mu_*)$, are represented
by formal pointed analytic homological superspaces, $(M_{\cB},0,\p)$ and,
respectively, $(M_{\fg},0, Q)$. Then the set
$$
S=\left\{\Gamma\in (\fg\ot m_{\cB})_{\tlo} \mid 
 \vec{\p}\Gamma =
\sum_{n=1}^{\infty} \frac{(-1)^{n(n+1)/2}}{n!}
\mu_n(\Gamma \bullet \ldots \bullet \Gamma)\right\}
$$
is just the set of all formal maps of pointed supermanifolds,
$$
\Gamma: (M_{\cB}, 0) \lon (M_{\fg},0),
$$
satisfying the equivariency condition
$$
d\Gamma (\p) = \Gamma^*(Q).
$$
The latter is precisely the ($L_{\infty}$-generalization of) the Master 
equation.

\sip

Both superpaces,  $(M_{\cB},0,\p)$ and  $(M_{\fg},0, Q)$, are foliated by 
integrable distributions,
\Beqrn
{\cal D}_{\p}&:=&\{X\in TM_{\cB}\mid X=[\p,Y]\ {\mathrm for\ some}\ Y\in 
TM_{\cB}\}, \\
{\cal D}_{Q}&:=&\{X'\in TM_{\fg}\mid X'=[Q,Y']\ {\mathrm for\ some}\ Y'\in 
TM_{\fg}\},
\Eeqrn
and $d\Gamma({\cal D}_{\p})\subset \Gamma^*({\cal D}_Q)$ for any $\Gamma\in S$.
Hence any such $\Gamma$ defines a map, $\hat{\Gamma}$, through 
the following commutative diagram,
$$
\diagram
M_{\cB}\dto \rto^{\Gamma} & M_{\fg}\dto \\
M_{\cB}/{\cal D}_{\p} \rto^{\hat{\Gamma}} & M_{\fg}/{\cal D}_{Q}\\
\enddiagram
$$
We say that two elements in $S$ are {\em gauge equivalent}, 
$\Gamma_1\sim \Gamma_2$, if $\hat{\Gamma}_1=\hat{\Gamma}_2$.
Infinitesimally, the gauge equivalence is given by
$$
\Gamma \sim \Gamma + (d+\vec{\p})g - 
\sum_{n=2}^{\infty} \frac{(-1)^{n(n+1)/2}}{(n-1)!}
\mu_n(g\bullet \Gamma \bullet \ldots \bullet \Gamma), \ \
\forall g\in (\fg\ot m_{\cB})_{\tln}.
$$

\bip

\paragraph{\bf 2.5.4. Remark.} If, for a derivation $\p:\cB\rar \cB$,
an element
 $\Gamma\in (\fg\ot m_{\cB})_{\tlo}$ satisfies the Master equation,
$$
\vec{\p}\Gamma =
\sum_{n=1}^{\infty} \frac{(-1)^{n(n+1)/2}}{n!}
\mu_n(\Gamma \bullet \ldots \bullet \Gamma),
$$
then, using the higher Jacobi identities (as in 2.5.1),
one gets an implication,
$$
\vec{\p}^{\,2}\Gamma = 0.
$$

\bip

\paragraph{\bf 2.5.5. Basic Theorem of 
$\Df$ormation Theory.} 
{\em Let $(\fg_1,Q_1)$ and $(\fg_2,Q_2)$ be two $L_{\infty}$-algebras
(in partucilar, dLie-algebras).
An $L_{\infty}$-morphism ${F}=\{{F}_n\}:\fg_1\rar\fg_2$
defines a natural transformation of the functors,
$$
\Ba{rccc}
{F}_*:& {\Df}^{\, *}_{\fg_1} &\lon& {\Df}^{\, *}_{\fg_2}\vspace{3mm} \\
& \Gamma &\lon & {F}_*(\Gamma):=\sum_{n=1}^{\infty} 
\frac{1}{n!}{F}_n(\Gamma,\ldots,\Gamma).
\Ea
$$

Moreover, if $F$ is a quasi-isomorphism, then $F_*$
is an isomorphism.
}

\sip

\Proof  Assume $\Gamma\in (\fg_1\ot m_{\cB})_{\tlo}$
satisfies the Master equation,
$$
d\Gamma (\p) = \Gamma^*(Q_1).
$$
The $L_{\infty}$-morphism $F$, when viewed as map,
$M_{\fg_1}\rar M_{\fg_2}$, of pointed formal manifold,
 satisfies,
$$
dF (Q_1) = F^*(Q_2).
$$
Hence $F_*(\Gamma)$, which the same as $F\circ \Gamma$,
obviously satisfies the Master equation in $(\fg_2,Q_2)$.
\hfill $\Box$

\bip

\paragraph{\bf 2.5.6. Smoothness Theorem.} {\em
The deformation functor $\Df^{ *}$ is unobstructed, i.e.\ 
for any dLie-algebra $\fg$ with finite-dimensional cohomology
$\bH(\fg)$,
the functor $\Df^{\, *}_{\fg}$ is versally representable by
a smooth pointed formal  $\dim \bH(\fg)$-dimensional homological
supermanifold $(\M_{\fg}, \p)$.

\sip

Moreover, the diffeomorphism class of $(\M_{\fg}, \p)$
is an invariant of $\fg$.
}

\sip

\Proof It is enough to show  that 
\Bi
\item[(i)]
there exists a versal element
$\Gamma\in k[[t]]\ot \fg$ and a differential $\p:k[[t]]\rar k[[t]]$
satisfying the Master equation
\Beq \label{master}
 d\Gamma+ \vec{\p}\Gamma + \frac{1}{2}[\Gamma\bullet \Gamma]=0,
\Eeq
\\
\item[(ii)]
the differential $\p$, when viewed as a vector field on the cohomological
supermanifold $\M_{\fg}$, is an invariant of $\fg$.
\Ei

\sip

Unless $\fg$ is formal, there is no quasi-isomorphism from $(\fg, 
[\ \bullet\ ],d)$ to its cohomology $(\bH(\fg), [\ \bullet\ ]_{\mathrm ind},
0)$. However, there {\em always}\,  exists a 
 $L_{\infty}$-structure, $\{\mu_*, {\mathrm with} \ \mu_{1}=0\}$, 
on $\bH(\fg)$ which is quasi-isomorphic, via some 
$L_{\infty}$-morphism $F$,  to $(\fg, 
[\ \bullet\ ],d)$.
Moreover,
this structure is defined uniquely up to
a homotopy.

\sip

Setting $\Gamma_{[1]}=\sum_i t^i [e_i]$, in the notations
of Sect.\ 2.2, we  define a derivation,
$\p:k[[t]]\rar k[[t]]$, by the formula
$$
\vec{\p}\Gamma_{[1]} =
\sum_{n=2}^{\infty} \frac{(-1)^{n(n+1)/2}}{n!}
\mu_n(\Gamma_{[1]} \bullet \ldots \bullet \Gamma_{[1]}).
$$
By Remark~2.5.4, this derivation satisfies $\p^2=0$. 
Hence $(\Gamma_{[1]}, \p)$
is a versal solution of the Master equation in $(\bH(\fg), \mu_*)$,
while $(F_*(\Gamma_{[1]}), \p)$ is, by Theorem~2.5.5,  a versal 
solution of the 
Master equation in $(\fg, 
[\ \bullet\ ],d)$. This proves claim (i).

\sip

The $L_{\infty}$-structure $\{\mu_*\}$ on $\bH\simeq \R^{p|q}$ 
is well-defined only up to 
a homotopy, $\{\eta_{(n)}: \odot^n \bH\rar \bH, n\geq 2,  
\widetilde{\eta_n}=\tlo\}$. Is is easy to check that  a homotopy change
of the induced $L_{\infty}$-structure
$$
\mu_* \stackrel{\eta_{(*)}}{\lon} \mu_*',
$$
does change the differential,
$$
\p\, {\lon}\, \p',
$$
 but in a remarkably geometric way,
$$
\p'= d\eta (\p),
$$
where $\eta:\R^{p|q}\rar \R^{p|q}$ 
is just a formal change of coordinates
$$
t^i \rar t'^i= t^i +   \sum_{j,k}\pm \eta^{\, i}_{(2)\,jk} t^jt^k + 
\sum_{j,k,l}\pm \eta^{\, i}_{(3)\,jkl} t^jt^kt^l + \ldots \ .
$$
Put another way, a homotopy change of $\mu_*$ affects only the coordinate
representation of the vector field $\p$ on  $\M_{\fg}$. As a geometric
entity, this is an invariant of $\fg$. \hfill $\Box$

\bip

\paragraph{\bf 2.5.7. Corollary.} {\em The derived category
of $L_{\infty}$-algebras with  finite-dimensional cohomology
is canonically equivalent to the (purely geometric) category
of  pointed formal homological supermanifolds, $(\M,\p,0)$,
with $\p$ satisfying $\p I\subset I^2$, $I$ being the ideal of 
the distinguished point $0\in \M$.}

\sip

\Proof It is well-known that each quasi-isomorphism 
of $L_{\infty}$-algebras is a homotopy equivalence. 
Thus the derived category of $L_{\infty}$-algebras
is equivalent to the their homotopy category. Then the required statement
follows immediately from an observation made in the proof of Theorem~2.5.6
that, for any homotopy class of $L_{\infty}$-algebras $[\fg]$,
the associated homotopy class of $L_{\infty}$-structures induced
on the cohomology $\bH(\fg)$ is isomorphically mapped into one and the same
homological manifold $(\M_{\fg},\p,0)$.
\hfill $\Box$

\bip

\paragraph{\bf 2.5.8. Remarks.} (i) 
The origin of the vector field $\p$ in 
Theorem~2.5.6  can be traced back to  Chen's 
power series connection \cite{Chen}. This will be made apparent
in Section~4 where we give another, perturbative, proof of 
the above Theorem. 
From now on we call $\p$ the {\em Chen's differential}\, or
{\em Chen's vector field}.

\sip

(ii) The higher order tensors $\mu_*$ induced on the cohomology
$\bH(\fg)$ by a $L_{\infty}$-quasi-isomorphism from
a dLie-algebra $\fg$ coincide precisely with the  Massey 
products when they are well-defined and univalued.
Thus the Chen's differential gives a compact (and invariant)
representation of the homotopy class of Massey products.

\bip

\paragraph{\bf 2.5.9. Extended Kuranishi moduli space.} 
Since the Chen's
vector field $\p$   on $\M$
is homological, the distribution
$$
{\cal D}_{\p}= \left\{ \X\in \T\, \mid \, \X\, =[\p, Y\,]\ \mbox{for some}\ 
Y\in \T\right\}
$$
is integrable (cf.\ Sect.\ 2.5.3). Indeed, the Jacobi identities imply
$$
\left[[\p,\X\,], [\p,\Y\,]\right]=\left[\p, \left[\X\,,[\p,
\Y\,]\right]\right].
$$
Consider an affine subscheme, 
$$
``{\mathsf Zeros}(\p)\mbox{''} := {\mathsf Spec}\  
k[[t]]/ <\p t^1, \ldots, \p t^{p+q}\,>,
$$
of zeros of the vector field $\p$. The distribution ${\cal D}_{\p}$ 
is tangent to $``{\mathsf Zeros}(\p)\mbox{''}$ since 
$[\p,[\p,\X\,]]=0$. We define 
the {\em extended Kuranishi space}, ${\mathcal K}_{\fg}$, 
as the so called non-linear homology \cite{BaKo,Ma,Ba} 
of the Chen differential, 
i.e.\ as the quotient
$``{\mathsf Zeros}(\p)\mbox{''}/{\cal D}_{\p}$. 
(For our purposes it is enough
to  understand the latter as ${\mathsf Spec}\,  k[[t]] \cap \Ker \p / 
<\p t^1, \ldots, \p t^{p+q} >$.)
This passage from $(\M,\p)$ to ${\mathcal K}_{\fg}$ 
establishes a clear link between
 the unobstructed deformation functor $\Df^{\, *}_{\fg}$
and the classical one $\Def^{\, *}_{\fg}$.

\sip

Kuranishi spaces, ${\mathsf K}_{\fg}$, originally 
emerged \cite{Ku2,GM} in the context 
of the deformation functor ${\mathrm Def}^0_{\fg}$ 
associated to a cohomologically
split $\Z$-graded dLie algebra $\fg$.
 It is 
not hard to see that ${\mathsf K}_{\fg}$ is a 
proper subspace of the extended Kuranishi space ${\mathcal K}_{\fg}$. 
We will see below  that for a rich class of dLie algebras $\fg$
--- the so-called (homotopy) Gerstenhaber algebras ---
the tangent sheaves to the smooth parts, 
${\mathcal K}_{\mathrm smooth}$, of the associated
extended Kuranishi spaces are 
canonically sheaves of associative algebras. This  
structure is { not}\, visible if working in the category
of original (non-extended) Kuranishi spaces ${\mathsf K}$ only.

\bip

\sip

\paragraph{\bf 2.6. Cohomological splitting.}
It is very easy to compute Chen's differential once 
a cohomological decomposition of a dLie-algebra $\fg$ under investigation
is chosen. The latter means the data $(i,p,Q)$, where
$i: \bH(\fg) \lon \fg$ is a linear injection,
$p:\fg \lon \bH(\fg)$  a linear surjection, and $Q:\fg \lon \fg$ an
 odd linear operator, all satisfying the conditions,
$$
p\circ i= \mbox{Id} = i\circ p \oplus dQ \oplus Qd, 
$$
in $\mbox{End}_k(\fg)$.

\sip

Such a decomposition of $\fg$ often occurs in (complex)
differential geometry \cite{Kod,Ku2}, where typical 
dLie-algebras come equipped with a norm $||\ ||$
and their cohomologies $\bH$ get identified 
with harmonic subspaces,  $\mbox{Harm}:=\Ker d\cap \Ker d^*\subset \fg$,
$d^*$ being the adjoint of $d$ with respect to $||\ ||$. The operator $Q$
is then $Gd^*$, where $G$ is 
the Green function of the Laplacian $\Box=d d^* + d^* d$.
In this situation the formal power series solution, $\Gamma$,  of 
the Master equation 
as well as the associated Chen's vector field
can be chosen to be 
convergent inducing, thereby, the structure of analytic (rather than formal)
homological supermanifold on $\M$.

\sip

It is not hard to check that, given a cohomological splitting of $\fg$,
the pair, $(\Gamma, \p)$,
 given recursively by
\Beqr
\Gamma_{[1]} &=& \sum_{i} t^{i}e_i 
\in \Ker Q\cap \Ker d \nonumber \\
\Gamma_{[2]} &=& - \frac{1}{2} Q [\Gamma_{[1]}(t)\bullet 
\Gamma_{[1]}(t)],\nonumber\\
\Gamma_{[3]} &=&  - \frac{1}{2} Q \left([\Gamma_{[1]}(t)\bullet \Gamma_{[2]}(t)] +
[\Gamma_{[2]}(t)\bullet \Gamma_{[1]}(t)]\right),\nonumber\\
\ldots && \nonumber\\
\Gamma_{[n]} &=& - \frac{1}{2} Q \left(\sum_{k=1}^{n-1}[\Gamma_{[k]}(t)\bullet
\Gamma_{[n-k]}(t)]
\right) \label{split} \\
\ldots && \nonumber
\Eeqr
and 
$$
\vec{\p}p(\Gamma_{[1]}):= - \frac{1}{2} p\left([\Gamma\bullet \Gamma]\right), 
$$
give an explicit versal solution of the Master equation in $\fg$.

\sip

The above power series for $\Gamma$ is well known in the 
classical Deformation Theory  
\cite{Kod,Ku2} where it plays a key role in constructing Kuranishi
analytic moduli spaces. This series is essentially an inversion
of the Kuranishi map \cite{Ku2} in the category of $L_{\infty}$-algebras.

\bip

\paragraph{\bf 2.7. Formality and flat structures.} 
If an algebra
$(\fg, d, [\ \bullet\ ])$ is formal, then,
as follows from the proof of
Theorem~2.5.6, the associated homological supermanifold $(\M_{\fg},\p,0)$
has a canonical flat structure $\nabla$. In the associated 
flat coordinates $t^i$, the Chen's vector field $\p$ has coefficients
which are polynomials in $t^i$ of order $\leq 2$. (This observation  can,
in fact, be made into a geometric criterion of formality.) More precisely,
the following is true.

\sip

{\bf 2.7.1. Theorem.} {\em
For any  formal dLie-algebra $\fg$ there is a canonical
isomorphism of sets},
$$
\left\{ \Ba{c} {\mathsf Flat\ structures\ on} \ \M_{\fg}\  
{\mathsf such\ that}\
 \nabla_X\nabla_Y\nabla_Z \p=0 \\
{\mathsf for\  any\ horizontal\
 vector\ fields}\ X,Y\ {\mathsf and}\ Z \Ea
\right\} \longleftrightarrow
\left\{ \Ba{c} {\mathsf Homotopy\  classes}\\
{\mathsf of \ formality\ maps}\Ea \right\}
$$ 

\bip

Let us choose a basis, $\{s_i, i=1, \ldots, p+q\}$, 
in the $(p|q)$-dimensional
vector superspace $\bH(\fg)$, and let $\{t^i\}$ be the associated
linear coordinates. We shall need, for a short time, a category
${\mathsf Artin}_{k[[t]]}$ consisting of Artin superalgebras
of the form 
$$
\cA_N:= {k[[t^1,\ldots,t^{p+q}]]}/< (t^1)^{N_1}\cdots 
(t^{p+q})^{N_{p+q}}>.
$$
Denoting the maximal ideal of such a superalgebra by $m_N$,
we set $(\fg\ot m_N)_{\mathsf versal}$
to be a linear subspace in $\fg\ot m_N$ consisting of even elements,
$\Gamma$, satisfying $\Gamma \bmod m_N=0$, $\Gamma \bmod m^2_N\in \Ker d$,
and
$$
\left(\Gamma \bmod m^2_N\right)\bmod \Img d = \sum_{i=1}^{p+q} t^i s_i.
$$
This set is invariant under the action of the gauge group
$\exp (\fg\ot m_N)_{\tln}$ (see Sect.\ 2.5).

\bip

{\bf 2.7.2. Lemma.} {\em For any formal dLie-algebra $\fg$ there is a canonical
isomorphism of sets,
$$
\lim_{\longleftarrow}    \frac{ \left\{ \Gamma\in (\fg\ot m_N)_{\mathsf versal}
\mid d\Gamma + \vec{\p} \Gamma + \frac{1}{2}[\Gamma\bullet \Gamma]=0\right\}}
{\mathsf gauge\ group}
= \frac{\mathsf Formality\ maps}{\mathsf homotopy\ equivalence}\, ,
$$
where the projective limit is taken over the category 
${\mathsf Artin}_{k[[t]]}$.
}

\bip

\Proof For any $\cA_N\in {\mathsf Artin}_{k[[t]]}$, 
the Master equation in the Lie algebra 
$(\bH(\fg)\ot \cA_N, [\ \bullet \ ]_{\mathrm ind})$ has a 
canonical versal solution,
$$
\Gamma_0 = \sum_{i=1}^{p+q} t^i s_i, \ \ \
\p= \sum_{i,j,k=1}^{p+q}(-1)^{\tl{j}(\tl{i}+1)}t^it^j C_{ij}^k \frac{\p}{\p
 t^k},
$$
where $C_{ij}^k$ are the structure constants of 
$[\ \bullet \ ]_{\mathrm ind}$.

\sip

If $\fg$ is formal, and $F=\{F_n: \odot^n \bH(\fg)\rar \fg, n=1,2,\ldots \}$ 
is a 
formality map, then, by Theorem~2.5.5, 
$$
\Gamma:= \sum_{n=1}^{\infty} \frac{1}{n!} F_n(\Gamma_0, \ldots, \Gamma_0),
\ \ \ {\mathsf the\ same\ \p}, 
$$
is a versal solution of the Master equation in $\fg\ot \cA_N$.

\sip

It is easy to check that an
arbitrary homotopy change, 
$$
F \rar {F}^h,
$$
 of the formality map, say the one induced by a  set of linear maps 
$h=\{h_n: \odot^n \bH(\fg)\rar \fg, \tl{h}_n=\tl{n},
 n=1,2,\ldots\}$, change the versal 
solution $\Gamma$ into a {\em gauge equivalent}\, one, 
$$
\Gamma \rar \Gamma^g,
$$
where
$$
g= \sum_{i}
h_1(e_i)t^i  +   
\sum_{i,j}\pm  h_2(e_i,e_j) t^it^j + 
\sum_{i,j,k}\pm  h_3(e_i,e_j,e_k) t^it^jt^k + \ldots \ .
$$

\sip

Hence there is a canonical map,
$$
\Ba{ccc}
\frac{\mathsf Formality\ maps}{\mathsf homotopy\ equivalence},
 & \lon & \frac{ \left\{ \Gamma\in (\fg\ot m_N)_{\mathsf versal}
\mid d\Gamma + \vec{\p} \Gamma + \frac{1}{2}[\Gamma\bullet \Gamma]=0\right\}}
{\mathsf gauge\ group} \vspace{4mm} \\
F=\{F_n: \odot^n \bH(\fg)\rar \fg, n=1,2,\ldots \}/\sim  & \lon &
\sum_{n=1}^{\infty} \frac{1}{n!} F_n(\Gamma_0, \ldots, \Gamma_0)/\sim
\Ea
$$
which implies (almost immediately) the desired result.
\hfill $\Box$

\bip

{\bf 2.7.3. Proof of Theorem~2.7.1.} Let ${\mathsf Diff}_0$ be the group
of all formal diffeomorphisms of $\M_{\fg}=\bH(\fg)$ into itself preserving
the origin, and set
$$
{\mathsf Diff}_{0,\p}:= \left\{\phi\in {\mathsf Diff}_0 \mid
\phi_*(\p)\ {\mathsf is\ quadratic\ in}\ t^i\right\}.
$$
Note that ${\mathsf Diff}_{0,\p}={\mathsf Diff}_0$ if the Chen's vector
field $\p$ vanishes. In general, 
$$
GL(p+q)\subseteq {\mathsf Diff}_{0,\p} \subseteq {\mathsf Diff}_{0}.
$$

\sip

There is an obvious isomorphism,
$$
\left\{ \Ba{c} {\mathsf Flat\ sructures\ on} \ \M_{\fg}\  
{\mathsf such\ that}\
 \nabla_X\nabla_Y\nabla_Z \p=0 \\
{\mathsf for\  any\ horizontal\
 vector\ fields}\ X,Y\ {\mathsf and}\ Z \Ea
\right\} =
\frac{{\mathsf Diff}_{0,\p}}{GL(p+q)}.
$$
On the other hand, by Theorem~4.2.2 (see below),
$$
\frac{{\mathsf Diff}_{0,\p}}{GL(p+q)}
= \lim_{\longleftarrow}    
\frac{ \left\{ \Gamma\in (\fg\ot m_N)_{\mathsf versal}
\mid d\Gamma + \vec{\p} \Gamma + \frac{1}{2}[\Gamma\bullet \Gamma]=0\right\}}
{\mathsf gauge\ group}.
$$
The final link in the chain of canonical isomorphisms
is provided by Lemma~2.7.2. 
\hfill $\Box$

\bip

{\bf 2.7.4. Corollary.} {\em  For any compact Calabi-Yau manifold,
there is a canonical isomorphism of sets,}
$$
\left\{ \Ba{c} {\mathsf Flat\ connections\ on\ Barannikov-Kontsevich's}\\
{\mathsf moduli\ space\ of \ extended\ complex\ structures}
\Ea \right\} \longleftrightarrow
\left\{ \Ba{c} {\mathsf Homotopy\  classes}\\
{\mathsf of \ formality\ maps}\Ea \right\}
$$

\bip



\section{Homotopy Gerstenhaber algebras}

\ \ \ \ \ \  \paragraph{\bf 3.1. Differential Gerstenhaber algebras. } 
A {\em differential Gerstenhaber algebra}, or shortly, a {\em dG-algebra},
is the data $(\fg,d, [\ \bullet\ ],\ \cdot\ )$ where
\Bi
\item[(i)] $(\fg,d, [\ \bullet\ ])$ is a $\Z$-graded dLie-algebra as
defined  in Sect.\ 2.2.1; 
\item[(ii)] $(\fg,d, \ \cdot\ )$ is a differential $\Z$-graded 
associative algebra with the product
$$
\begin{array}{rccc}
 \cdot: & \fg\ot \fg & \lon & \fg \\
& a \ot b & \lon &  a\cdot b,
\end{array}
$$
having degree $0$;
\item[(iii)] the binary operations $[\ \bullet\ ]$ and $\cdot$
satisfy the odd Poisson identity,
$$
[a\bullet (b\cdot c)]=[a\bullet b]\cdot c +(-1)^{(\tl{a}+1)\tl{b}}b\cdot
[a\bullet c],
$$
\Ei
for all homogeneous $a,b,c\in \fg$.

\sip

A dG-algebra is called 
{\em  graded commutative} if such is the dot product.
\sip

The {\em identity}\, in $\fg$ is an even element $e_0$ such
that $de_0=0$, $e_0\cdot a=a\cdot e_0=a$ and $[e_0\bullet a]=0$  for any $a\in \fg $. 
It defines
a cohomology class $[e_0]$ in $\bH$, and a constant vector field  
on $\M$ which we denote by $e$.



\bip

\paragraph{\bf 3.1.1. Remark.}
If $\fg$ is a unital dG-algebra, then a versal solution, $\Gamma$,
of the Master equation (\ref{master}) in $\fg$ can (and will)  be always 
normalized in such a way that
$$
\vec{e}\, \Gamma =e_0.
$$

\bip

\paragraph{\bf 3.1.2. Differential Gerstenhaber-Batalin-Vilkovisky algebras.} 
Let $(\fg, \ \cdot\ )$ be a 
$\Z$-graded commutative associative 
algebra over a field $k$. Let us say that the zero operator, $0:\fg\rar \fg$,
is of order $-1$, and let us denote the linear operator, $x\rar a\cdot x$,
 of left multiplication by an element $a\in \fg$ by $l_a$. A homogeneous 
linear operator, $D:\fg\rar \fg$, is said
to be an {\em operator of order $k$} if the operator
$[D,l_a]$ is of order $k-1$ for any homogeneous $a$ in $\fg$.

\sip

 Assume now that $(\fg, \ \cdot\ )$  comes equipped with 
\begin{itemize}
\item[(i)] a degree $+1$ linear
operator, $d:\fg\rar \fg$, of order 1, and \\
\item[(ii)]   a degree $-1$ linear  operator, 
$\Delta:\fg\rar \fg$, of order $2$,
\end{itemize}
satisfying the conditions,
$$
d^2=0, \ \ \ \ \ \Delta^2=0, \ \ \ \ \ d\Delta+\Delta d=0.
$$
In this case the data
$$
\left(\fg, d, \ \cdot\ , [\ \bullet\ ]\right)
$$
with
$$
[a\bullet b]:= (-1)^{\tl{a}}\Delta(a\cdot b) - (-1)^{\tl{a}}(\Delta a)\cdot b
- a\cdot (\Delta b), \ \ \ \ \ \forall a,b\in \fg,
$$
defines a dG-algebra \cite{Ma}. The dG-algebras arising in this way
are  often called {\em exact}\, or dGBV-algebras. 

\bip

\paragraph{\bf 3.1.3. Example (complex manifolds).} For any  $n$-dimensional 
comlex manifold $M$ the differential 
Lie algebra of Example~2.3.1,
$$
\fg=\left(\Gamma(M, \Lambda^{\bullet} T_M \ot \Lambda^{\bullet} 
\overline{T}_M^*), 
\ \bar{\p}\ ,\ [\ \bullet\ ]\ \right),
$$
equipped with a supercommutative product,
$$
\begin{array}{rccc} 
 \wedge : & \Gamma(M,\Lambda^{i_1} T_M \ot 
\Lambda^{j_1}\ovl{T}_M^*)
\, \times \, \Gamma(M,\Lambda^{i_2} T_M \ot \Lambda^{j_2}\ovl{T}_M^*) 
& \rar & \Gamma(M,
\Lambda^{i_1+i_2} T_M \ot \Lambda^{j_1+j_2}\ovl{T}_M^*)  \\
& \X_1\ot \ovl{w}_1 \, \times \X_2\ot \ovl{w}_2  & \rar & 
(-1)^{\tl{j}_1\tl{i}_2} (\X_1\wedge \X_2) \ot \ovl{w}_1\wedge \ovl{w}_2,
\end{array}
$$
is a unital $\Z$-graded commutative dG-algebra.

\sip

If $M$ admits a nowhere vanishing global holomorphic volume form, 
$\Omega\in \Gamma(M,
\Omega^{n}_M)$, then the above dG-algebra is actually exact 
\cite{Tian,To,BaKo} 
 with $\Delta$ being the composition,
$$
\Delta:\ \Lambda^i T_M \stackrel{i_{\Omega}}{\lon} 	\Omega^{n-i}_M 
\stackrel{\p}{\lon} \Omega^{n-i+1}_M  \stackrel{i_{\Omega}^{-1}}{\lon}
 \Lambda^{i-1} T_M.
$$
Here $i_{\Omega}: \Lambda^{\bullet} T_M \rar \Omega^{\bullet}_M$ 
is the natural isomorphism
given by contraction with the holomorphic volume form, and $\p$ is the 
$(1,0)$-part of the de Rham operator.

\bip

\paragraph{\bf 3.1.4. Example (symplectic manifolds).} 
For any symplectic manifold $(M,\omega)$ the 
dLie 
algebra of Example~2.3.5,
$$
\fg=\left(\Gamma(M, \Omega^{\bullet}_{\R}),\  d\ ,\  [\ \bullet\ ]_{\omega}\, 
\right),
$$
together with a graded commutative product, $a\cdot b:= a\wedge b$, is a 
unital $\Z$-graded dG-algebra.
Moreover, it is a dGBV-algebra with the 2-nd order differential given by
$$
\Delta|_{\Omega^k_{\R}} =(-1)^{k+1}*d*.
$$
Here  $*: \Omega^k_{\R} \rar\Omega^{2m-k}_{\R}$ is the 
symplectic analogue of the Hodge
duality operator
defined by the condition, $\be\wedge (*\al)= \langle \be, \al \rangle
\omega^m/m!$, with $\langle\, , \, \rangle$ being  the pairing between
$k$-forms induced by the symplectic form.

\bip

\paragraph{\bf 3.1.5. Example (vector bundles).} 
Let $M$ be a complex manifold, 
and $\pi: E\rar M$ a holomorphic vector bundle. 
There is a complex of $\Z$-graded sheaves 
$(\odot^{\bullet}E \ot  \Lambda^{\bullet}E^*,\Delta)$,
$$
\ldots \stackrel{\Delta}{\lon} \odot^{k+1} E\ot \Lambda^{l+1} E^* 
\stackrel{\Delta}{\lon}
\odot^k E\ot \Lambda^{l} E^* \stackrel{\Delta}{\lon}
\odot^{k-1} E\ot \Lambda^{l-1} E^* \stackrel{\Delta}{\lon} \ldots,
$$
where the differential $\Delta$ is just the contraction, 
the $\Z$-grading is induced
from the one on $\Lambda^{\bullet} E$,  and we set 
$\odot^k E=\Lambda^k E^* =0$ for $k< 0$. It is easy to see that 
$\Delta$ is a linear
operator of oder 2 with respect to the natural supercommutative product,
$$
\begin{array}{ccc} 
\odot^{\bullet}E \ot  \Lambda^{\bullet}E^*
\, \times \, \odot^{\bullet}E \ot  \Lambda^{\bullet}E^*
& \stackrel{\cdot}{\lon} & \odot^{\bullet}E \ot  \Lambda^{\bullet}E^* \\
 (a_1\ot b_1^*) \, \times (a_2\ot b_2^*)  & \lon & 
(a_1\odot a_2)\ot (b_1^*\wedge b_2^*).
\end{array}
$$

\sip

Hence the data
$$
\fg= \left(\bigoplus_k \fg^k, 
\fg^k:= \bigoplus_{i+j=k} 
\Gamma(M, 
\odot^{\bullet}E 
\ot \Lambda^{i}E^* \ot\overline{\Omega}^j_{M}),\ \Delta\ , \ 
\bar{\p}\right)
$$
is a unital dGBV-algebra.  It extends the dLie algebra of 
Example~2.3.6,
as the following calculation shows.

\bip

\paragraph{\bf Proposition.} {\em The bracket $[\ \bullet\ ]$ on 
$\odot^{\bullet}E \ot  \Lambda^{\bullet}E^*$, when restricted 
to $E\ot E^*$, coincides, up to a sign, 
with the usual commutator of morphisms.}

\sip

\Proof Let us consider a pair of germs, $C_1= a_1\ot b_1^*$ and 
$C_2=a_2\ot b_2^*$, in the same stalk of $E\ot E^*$. Then
\Beqrn
[C_1\bullet C_2] &=& (-1)^{\tl{C}_1} \Delta(C_1\cdot C_2) - (-1)^{\tl{C}_1}  
\Delta(C_1)\cdot C_2 - C_1\cdot \Delta(C_2) \\
&=& - \Delta\left((a_1\odot a_2)\ot (b_1^*\wedge b_2^*)\right)
+ \Delta(a_1\ot b_1^*)\cdot (a_2\ot b_2) 
- (a_1\ot b_1) \cdot \Delta(a_2\ot b_2^*) \\
&=& - \langle a_1, b_1^*\rangle \, a_2\ot b_2^* +  
\langle a_1, b_2^*\rangle\, a_2\ot b_1^*
-  \langle a_2, b_1^*\rangle\, a_1\ot b_2^* \\
&&
 +  \langle a_2, b_2^*\rangle\, a_1\ot b_1^*  +\, \langle a_1, b_1^*\rangle\, 
a_2\ot b_2^* -  \langle a_2, b_2^*\rangle\, a_1\ot b_1^*\\
&=& \langle a_1, b_2^*\rangle\, a_2\ot b_1^*
-  \langle a_2, b_1^*\rangle\, a_1\ot b_2^* \\
&=& -\left(C_1C_2 - C_2C_1\right),
\Eeqrn
where angular brackets stand for the usual pairing between a vector and a 1-form.
\hfill $\Box$

\bip

There is a  problem with the constructed  dGBV-algebra ---  
its cohomology 
may not be finite dimensional even for compact manifolds. It is can be resolved
by passing to its dGBV-subalgebra,
$$
\fg_E= \left(\bigoplus_k \fg^k, 
\fg^k:= \bigoplus_{i+j=k} 
\Gamma(M, 
\odot^i E \ot \Lambda^{i}E^* \ot\overline{\Omega}^j_{M}),\ \Delta\ , \ 
\bar{\p}\right).
$$

If necessary, the asymmetry of $E$ and $E^*$ can be eliminated by taking
the tensor product $\fg_E\ot \fg_{E^*}$. 

\bip

\paragraph{\bf 3.1.6. Example (Hochschild cohomology).} Let $A$
be an associative algebra over a field $k$. The $\Z$-graded vector space
of {\em Hochschild cochains},
$$
C^{\bullet}(A,A) := \bigoplus_{n=0}^{\infty} \Hom_k(A^{\ot n}, A),
$$
can be made into a {\em Hochschild complex}\, with the
differential, $d: C^{n}(A,A)\rar C^{n+1}(A,A)$, given by
\Beqrn
(df)(a_1\ot \ldots \ot a_{n+1}) &:=&
a_1 f(a_2\ot \ldots \ot a_{n+1})
+ \sum_{i=1}^n (-1)^I f(a_1\ot \ldots \ot a_i a_{i+1}\ot \ldots\ot a_{n+1})\\
&&  +\, (-1)^{n+1}f(a_1\ot \ldots \ot a_n)  a_{n+1}
\Eeqrn
for any $f\in C^n(A,A)$.

\sip

One can define two binary operations,
$C^{\bullet}(A,A) \ot C^{\bullet}(A,A) \rar C^{\bullet}(A,A)$,
the degree 0 {\em dot product},
$$
(f\cdot g)(a_1\ot \ldots \ot a_{k+l}):= (-1)^{kl}f(a_1\ot \ldots \ot a_{k})
g(a_1\ot \ldots \ot a_{l}), \ \ \forall f\in C^k(A,A), g\in C^l(A,A),
$$
and the degree $-1$ {\em  bracket},
$$
[f\bullet g]:= f\circ g - (-1)^{(k+1)(l+1)}g\circ f,
$$
where
$$
(f\circ g) (a_1\ot \ldots \ot a_{k+l-1}):= 
\sum_{i=1}^{k-1} (-1)^{(i+1)(l+1)} 
 f(a_1\ot \ldots \ot a_i\ot g(a_{i+1}\ot \ldots \ot a_{i+l})
\ot \ldots \ot a_{k+l-1}).
$$ 
These two make the Hochschild complex into a $\Z$-graded
differential  associative algebra and a differential (odd) Lie algebra
respectively. Though $(C^{\bullet}(A,A), d, [\ \bullet\ ],\ \cdot\ )$
is not a dG-algebra, it is a remarkabale fact that 
the  
associated Hochschild cohomology,
$$
{\mathrm Hoch}^{\bullet}(A,A)= \frac{\Ker d}{\Img d},
$$
carries the structure of graded commutative dG-algebra with respect to the
naturally indiced dot product, Lie bracket and the zero differential.

\bip

\paragraph{\bf 3.2. $A_{\infty}$-algebras.}
A {\em strong homotopy algebra}, or shortly $A_{\infty}$-{\em algebra},
is by definition a vector superspace $V$ equipped with linear maps,
$$
\Ba{rcccc}
\mu_k: & \ot^k V & \lon & V & \\
& v_1\ot\ldots \ot v_k & \lon & \mu_k(v_1, \ldots, v_k), & \ \ \ \ k\geq 1,
\Ea
$$
of parity $\tilde{k}$ satisfying, for any $n\geq 1$ and
 any $v_1, \ldots, v_n \in V$,
 the following
{\em higher order associativity conditions},
\Beq \label{id}
\sum_{k+l=n+1}\sum_{j=0}^{ k-1} (-1)^r \mu_k\left(v_1,\ldots,v_j,
\mu_l(v_{j+1},\ldots, v_{j+l}), v_{j+l+1}, \ldots, v_n\right)=0,
\Eeq
where $r=\tilde{l}(\tv_1 +\ldots + \tv_j) + \tj(\tll-1) + (\tk-1)\tll$
and $\tv$ denotes the parity of $v\in V$.

\vse

Denoting $dv_1:=\mu_1(v_1)$ and $v_1\cdot v_2 := \mu_2(v_1,v_2)$, we can spell
the first three conditions from the above infinite series as follows,
\begin{description}
\item[$n=1$:]\ \  $d^2=0$,\\
\item[$n=2$:]\ \  $d(v_1\cdot v_2)= (dv_1)\cdot v_2 + (-1)^{\tv_1} v_1\cdot
(dv_2)$,\\
\item[$n=3$:]\ \ $v_1\cdot (v_2\cdot v_3) - (v_1\cdot v_2)\cdot v_3 =
d \mu_3(v_1,v_2,v_3) + \mu_3(dv_1, v_2, v_3) + (-1)^{\tv_1} \mu_3(v_1, dv_2,
v_3) + \mbox{\hspace{5mm}$(-1)^{\tv_1+ \tv_2}\mu_3(v_1,v_2,dv_3)$}$,\\
\end{description}

Therefore $A_{\infty}$-algebras with $\mu_k=0$ for $k\geq 3$
are nothing but the differential associative superalgebras with the
differential $\mu_1$ and the associative multiplication 
$\mu_2$. If, furthermore, $\mu_1=0$, one recovers  the usual associative
superalgebras.

\sip

There is a (finer) $\Z$-graded version of the above definition in which 
the maps $\mu_n$ are required to be homogeneous (usually of degree $n-2$) 
with respect to the given $\Z$-grading on $V$.

\bip

\paragraph{\bf 3.2.1. Identity.}
An element $e$ in the $A_{\infty}$-algebra is called the {\em identity}\, if
$\mu_1(e)=0$, $\mu_2(e,v)=\mu_2(v,e)=v$ and $\mu_n(v_1,\ldots, e, \ldots, 
v_{n-1})=0$ for all $n\geq 3$ and arbitrary $v, v_1, \ldots, v_{n-1}\in V$.

\bip

\paragraph{\bf 3.2.2. Homotopy classes of $A_{\infty}$-algebras.}
For a pair of $A_{\infty}$-algebras,  $(V, \mu_*)$ and $(\tilde{V}, 
\tilde{\mu}_*)$,
 there is a natural
notion of a $A_{\infty}$-morphism from $V$ to $\tilde{V}$ which is, 
by definition,
a set of linear maps 
$$
F=\{f_n:V^{\otimes n}\lon \tilde{V}, \ n\geq 1 \},
$$
of parity $\tl{n}+1$ (or of degree  $1-n$ in the $\Z$-graded case)
which satisfy
$$
\sum_{1\le k_1<k_2<\ldots<k_i=n}(-1)^{i+r}
\tmu_i(f_{k_1}(v_1,\ldots,v_{k_1}),f_{k_2-k_1}(v_{k_1+1},\ldots,v_{k_2}),
\ldots,f_{n-k_{i-1}}(v_{k_{i-1}+1},\ldots,v_n))
$$
$$
=\sum_{k+l=n+1}\sum_{j=0}^{k-1}(-1)^{l(\tv_1+\ldots+\tv_j+n)+j(l-1)}
f_k(v_1,\ldots,v_j,\mu_l(v_{j+1},\ldots,v_{j+l}),v_{j+l+1},\ldots,v_n).
$$
The first three floors in the above infinite tower are
\begin{description}
\item[$n=1$:]\ \  $\tmu_1=\mu_2=:d$,\\
\item[$n=2$:]\ \  $\tmu_2(v_1, v_2)= \mu_2(v_1,v_2) + (df_2)(v_1,v_2)$, \\
\item[$n=3$:]\ \ $\tmu_3(v_1,v_2,v_3) + \tmu_2(f_2(v_1,v_2),v_3)
- (-1)^{\tv_1}\tmu_2(a_1, f_2(v_3,v_4))=$\\
\hspace{42mm} $\mu_3(v_1,v_2,v_3) - f_2(\mu_2(v_1,v_2),v_3)
+ f_2(v_1, \mu_2(v_3,v_4)) + (df_3)(v_1,v_2,v_3),$\\
\end{description}
where we naturally extended the differential $d:V\rar V$ to $d: \ot^k V^*\ot V 
\rar \ot^kV^*\ot V$ (so that, for example,
$(df_2)(v_1,v_2)= df_2(v_1,v_2) + f_2(dv_1,v_2) + (-1)^{\tv_1} f_2(v_1,dv_2)$)

\sip

A morphism $F=\{f_n\}$ of the $A_{\infty}$-algebra $(V,\mu_*)$
to itself is called a {\em homotopy}\,
if $f_1$ is an isomorphism. 
If $(V,\mu_*)$ has the identity $e$, then by a homotopy
of $(V,\mu_*,e)$ we understand a homotopy of $(V,\mu_*)$ satisfying 
the additional conditions, $f_n(v_1, \ldots, e, \ldots, v_{n-1})=0$ for all $n\geq 2$
and arbitrary $v_1,\ldots,v_{n-1}\in V$.

\sip 

It is not hard to see that
homotopy defines an equivalence relation in the set of 
all possible (unital) $A_{\infty}$-structures on a given vector superspace $V$. 

\bip

\paragraph{\bf 3.2.3. Remark.}
For future reference we rewrite the $n$-th order associativity condition
(\ref{id}) as
$$
\Lambda_n(v_1,\ldots,v_n) = (d\mu_n)(v_1,\ldots,v_n)
$$
where
$$
\Lambda_n(v_1,\ldots,v_n) :=\sum_{k+l=n-1\atop k,l\geq 1}
(-1)^{r'}\mu_{k+1}(v_1,\dots,v_j,
\mu_{l+1}(v_{j+1}, \ldots, v_{j+l+1}), v_{j+l+2}, \ldots, v_n)
$$
and $r'=(l+1)(\tv_1+\ldots+\tv_j) + jl + k(l+1) +1$.

\bip 

\paragraph{\bf 3.2.4. Remark.} It follows from (\ref{id}) for $n=3$ that
 the cohomology,
$$
H(V):=\frac{\Ker\, \mu_1}{\Img\, \mu_1}
$$
of a (unital) $A_{\infty}$-algebra $(V,\mu_*)$ is 
canonically a (unital) {\em associative}\, algebra. 
Moreover, a homotopy class of (unital) $A_{\infty}$-structures on $V$ 
induces one and the same
structure of (unital) associative algebra on $H(V)$.

\bip

\paragraph{\bf 3.2.5. The bar construction.}
There is a  conceptually better interpretation \cite{St1} of an 
$A_{\infty}$-structure on the vector superspace $V$ as a co-differential 
on the bar-construction of $V$. Here are the details: 
\begin{itemize}
\item[(i)] The vector space 
$$
{\mathsf B}(V):= \bigoplus_{n=1}^{\infty} \left(V[1]\right)^{\ot n}
$$
is naturally a co-algebra with the co-product given by
$$
\Delta(w_1\ot \ldots \ot w_n)= \sum_{i=1}^n \left(w_1\ot \ldots \ot w_i\right)
\ot \left(w_i\ot \ldots \ot w_n\right).
$$
\item[(ii)] A linear map $Q: {\mathsf B}(V) \rar {\mathsf B}(V)$
is said to be a {\em co-derivation}\, if $\Delta Q= Q\ot \Id
+ \Id\ot Q$. There is a one-to-one correspondence between such co-derivations
and Hochschild cochains understood as elements of 
$\Hom({\mathsf B}(V), V)$.

\item[(iii)] A homogeneous (of degree $-2$) 
Hochschild cochain $\mu_*: {\mathsf B}(V)\rar V$ 
defines an $A_{\infty}$-structure on $V$ if and only if the associated
co-derivation $Q$ is a {\em co-differential}, i.e.\ satisfies $Q^2=0$.

\end{itemize}

In this setup, a morphism 
$(V, \mu_*)\rar (\tilde{V}, \tilde{\mu}_*)$ as in 3.8.1 is precisely
 a morphism of the associated bar-constructions respecting 
co-differentials.

\bip

\sip

\paragraph{\bf 3.3. $C_{\infty}$-algebras.}
This notion is a supercommutative analogue of the notion of 
$A_{\infty}$-algebra.

\sip

Let $V$ be a $\Z$-graded vector space and ${\mathsf B}(V)$ 
its bar construction. One can make the latter into an associative and graded 
commutative  algebra by defining the {\em shuffle tensor product},
$\circledast:{\mathsf B}(V)\ot {\mathsf B}(V){\rar} {\mathsf B}(V)$,
as follows
$$
(w_1\ot \ldots \ot w_k)\circledast(w_{k+1}\ot \ldots \ot w_n)
:= \sum_{\sigma\in Sh(k,n)}e(\sigma; w_1,\ldots, w_n)
w_{\sigma(1)}\ot\ldots \ot w_{\sigma(n)}.
$$
Here we used the notations explained in Sect.\ 2.4.

\sip

By definition \cite{GJ},
a {\em strong homotopy commutative algebra}, or shortly, 
$C_{\infty}$-algebra is an $A_{\infty}$-algebra $(V, \mu_*)$
such that the associated 
Hochschild cochain $\mu_*: {\mathsf B}(V)\rar V$ factors through
the composition\footnote{Such cochains are often called {\em Harrison}\,
cochains.}
$$
\diagram
\mu_*: \ {\mathsf B}(V) \ \ \ 
\rto^{\hspace{-8mm}{\mathrm natural \atop projection}} & 
 \ \ \ {\mathsf B}(V)/{\mathsf B}(V)\circledast{\mathsf B}(V) \rto & V.  
\enddiagram
$$
This implies, in particular, that 
$$
\mu_2(v_1,v_2)=(-1)^{\tv_1\tv_2}\mu_2(v_2,v_1)
$$
for any $v_1,v_2\in V$.

\sip

One defines notions of unital $C_{\infty}$-algebras, of a morphism of
$C_{\infty}$-algebras, of their homotopy etc.\ in the same way as 
in the $A_{\infty}$-case.

\bip

\sip

\paragraph{\bf 3.4. $G_{\infty}$-algebras.}
Let $V$ be a $\Z$-graded vector space and let
$$
{\mathsf Lie}(V[1]^*)= \sum_{k=1}^{\infty}{\mathsf Lie}^k(V[1]^*)
$$
the free graded Lie algebra generated by the shifted dual vector space
 $V[1]^*$, i.e.
$$
{\mathsf Lie}^1(V[1]^*):= V[1]^*,  \ \ \  
{\mathsf Lie}^k(V[1]^*):= \left[ V[1]^*, {\mathsf Lie}^{k-1}(V[1]^*)
\right].
$$
The Lie bracket on ${\mathsf Lie}(V[1]^*)$ extends in a usual way
to the skew-symmetric associative algebra
$$
\wedge^{\bullet} {\mathsf Lie}(V[1]^*) = \sum_{k=0}^{\infty} \wedge^k
{\mathsf Lie}(V[1]^*),
$$
making the latter into a Gerstenhaber algebra.

\bip

\paragraph{\bf 3.4.1. Definition \cite{Ta1,TT}.}
{\em A homotopy Gerstenhaber algebra}, or shortly 
{\em $G_{\infty}$-algebra} is a graded vector space $V$ together with
a degree one linear operator
$$
Q:  \wedge^{\bullet} 
{\mathsf Lie}(V[1]^*) \lon \wedge^{\bullet} {\mathsf Lie}(V[1]^*)
$$
such that $Q^2=0$ and $Q$ is a derivation with respect to both the product
and the bracket.

\sip

A $G_{\infty}$-{\em morphism}, $V\rar V'$, of $G_{\infty}$-algebras
is by definition a morphism,
$(\wedge^{\bullet} {\mathsf Lie}(V[1]^*), Q)\rar
(\wedge^{\bullet} {\mathsf Lie}(V'[1]^*),Q')$, of associated
differential Gerstenhaber algebras.

\sip

The definition 3.4.1 makes sense only in the case when $V$  is 
finite-dimensional.
However, an obvious dualization fixes the problem \cite{TT}:
\Bi
\item[(i)] The dual of ${\mathsf Lie}(V[1]^*)$ can be
identified with the quotient 
 ${\mathsf B}(V)/ {\mathsf B}(V)\circledast {\mathsf B}(V)$,
$\circledast$ being the shuffle tensor product. \\
\item[(ii)] Derivations of $\wedge^* {\mathsf Lie}(V[1]^*)$
can be identified with arbitrary collections of linear maps,
$$
m^*_{k_1,\ldots, k_n}: V[1]^* \lon  {\mathsf Lie}^{k_1}(V[1]^*)
\wedge\ldots\wedge {\mathsf Lie}^{k_n}(V[1]^*),
$$
which upon dualization go into linear homogeneous maps,
$$
m_{k_1,\ldots, k_n}: \frac{V^{\ot k_1}}{\mathrm shuffle\ products}
 \odot \ldots \odot \frac{V^{\ot k_n}}{\mathrm shuffle\ products}
\lon V,
$$
of degree $3-n-k_1-\ldots-k_n$.
\item[(iii)] the condition $Q^2=0$ translates  into a well-defined
set of quadratic equations for $m_{k_1,\ldots, k_n}$ which say, 
in particular, that $m_1$ is a differential on $V$ and that the product,
$v_1\cdot v_2:=(-1)^{\tv_1}m_2(v_1,v_2)$, together with 
the Lie bracket,
 $[v_1\bullet v_2]:=-(-1)^{\tv_1}
m_{1,1}(v_1,v_2)$, satisfy the Poisson identity up to a homotopy
given by $m_{2,1}$. Hence the associated
cohomology space $\bH$ is a graded commutative Gerstenhaber algebra
with respect to the  binary operations induced by  $m_2$ and $m_{1,1}$.
\Ei

The {\em identity}\, in a $G_{\infty}$-algebra $V$ is an even
element $e$ such that all $m_{k_1, \ldots, k_n}(\ldots,e,\ldots)$
vanish except $m_2(e,v)=v$.

\bip

\paragraph{\bf 3.4.2. Theorem-construction.} {\em There 
is a canonical functor from the derived category of 
unital $G_{\infty}$-algebras with finite-dimensional cohomology 
to the category of $F_{\infty}$-manifolds.}

\sip

\Proof Since each quasi-isomorphism of $G_{\infty}$-algebras is an
equivalence relation, the derived category of $G_{\infty}$-algebras
 coincides with their homotopy category.

\sip

We construct the desired functor,
$$
\diagram
{\Ba{c} \mathsf{Derived\ category}\\
{\mathsf of}\ 
G_{\infty}\ {\mathsf algebras}\Ea} \framed & \stackrel{F_{\infty}}
{\lon}
&
{\Ba{c} F_{\infty}\\  {\mathsf manifolds }\Ea} \framed 
\enddiagram
$$
in two steps.

\sip

{\em Step 1.}
Suppose we are given a 
homotopy class, $[\ ]$, of $G_{\infty}$-structures on a graded
vector space $V$. By Kontsevich's Lemma 1 in \cite{Ko3}, 
a cohomological splitting of the complex $(V, m_1)$
transfers $[\ ]$ into a 
homotopy class, $[ \wedge^{\bullet} 
{\mathsf Lie}(\bH[1]^*), Q]$, of minimal $G_{\infty}$-algebras
on the finite-dimensional cohomology space of the above complex.
Moreover, this class does not depend on the choice of a particular
cohomological splitting, and it is homotopy equivalent to the 
original one.

\sip

{\em Step 2.} Let ${\cal I}$ be the multiplicative ideal in  
$\wedge^{\bullet} {\mathsf Lie}(\bH[1]^*)$  generated by the commutant
of  ${\mathsf Lie}(\bH[1]^*)$. Any differential $Q$ from the induced
homotopy class preserves this ideal and induces, through the quotient
$\wedge^{\bullet} {\mathsf Lie}(\bH[1]^*)/{\cal I}$,
a homotopy class of  $L_{\infty}$-structures on $V$ which, by 
Corollary~2.5.7, can be identified with an odd 
vector field $\p$ on the associated cohomological supermanifold $\M$
satisfying $\p I\subset I^2$ and $[E,\p]=\p$, $I$ being the ideal
of the distinguished point $0\in \M$ and $E$ the Euler vector field.
We claim that the rest of the data listed  in Definition~1.1 gets induced
on $\M$ through the quotient
$\wedge^{\bullet} {\mathsf Lie}(\bH[1]^*)/{\cal I}^2$. Indeed,
what is left of a differential $Q$ on this quotient can be described
as a collection of tensors, $m_{k,1,\ldots,1}$, which, in a basis
$\{e_a\}$ of $\bH$, are represented by their components,
$\mu^a_{b_1\ldots b_k, c_1,\ldots, c_l}$, $k\geq 1, l\geq 0$. 
The Chen's vector field
$\p$ 
and the tensors $\mu_k$ defining the structure of a $C_{\infty}$-algebra
on the tangent sheaf $\T$ are then given by  formal power series,
$$
\p= \sum_{l\geq 0}\pm\mu^a_{b_1,c_1,\ldots,c_l}t^{b_1}t^{c_1}\ldots
t^{c_l}\frac{\p}{\p t^a}
$$
and
$$
\mu^a_{b_1\ldots b_k}=  
\sum_{l\geq 0}\pm\mu^a_{b_1\ldots b_k,c_1,\ldots,c_l}t^{c_1}\ldots
t^{c_l}.
$$
where $t^a$ are the associated linear coordinates on $\M$ to which we
assign degree $2-|e_a|$.
It is easy to see that the $G_{\infty}$-identities for $m_{k,1,\ldots,1}$
get transformed into  the right identities for the tensor fields
$\p$ and $\mu_k$ on $\M$. 
This completes the construction. \hfill $\Box$

\sip

\bip

\paragraph{\bf 3.4.3. Corollary.} {\em 
For any unital $G_{\infty}$-algebra with finite dimensional
cohomology,
the tangent sheaf
 to the smooth part of the extended
Kuranishi space, $\cM=``{\mathsf zeros}({\p})\mbox{''}/ 
\Img\mu_1$, 
 is canonically a sheaf of
induced (unital) associative algebras}.

\bip

It will be interesting to find out when $\cM_{\mathrm smooth}$
with its canonically induced structure 3.4.3
is an $F$-manifold in the sense of Hertling and Manin \cite{HM}.

\bip
\sip

\paragraph{\bf 3.5. Remark.} Different ``resolutions'' of the 
chain operad in the little disk operad give different notions
of homotopy Gerstenhaber algebra \cite{V}.
 The definition~3.4.1
is the most canonical one. However, the functor $F_{\infty}$ is not an 
equivalence in this case.

\sip

The proof of Theorem~3.4.2 suggests one more version:
{\em a reduced homotopy Gerstenhaber algebra}\, is a graded vector space
$V$ together with the structure of $G_{\infty}$-algebra
such that all composition maps $m_{k_1,\ldots, k_n}$ vanish
except $m_{k_1,1,\ldots, 1}$.
The derived category of such algebras is equivalent to the category
of $F_{\infty}$-manifolds (cf.\ Theorem~2.5.7).

\bip

\paragraph{\bf 3.6. Formality and Gauss-Manin connections.} 
A {\em pre}-${\mathit Frobenius}_{\infty}$
{\em manifold}\, is
the data $(\M, E,\nabla, \p, [\mu_{*}], e)$, where
\begin{itemize}
\item[(i)] $\M$ is a formal pointed $\Z$-graded manifold, \\
 \item[(ii)] $E$ is the Euler vector field on $\M$, 
$Ef:=\frac{1}{2}|f|f$, for all homogeneous functions on $\M$
of degree $|f|$,
\item[(iii)] $\nabla$ is a flat torsion-free affine connection,
called the Gauss-Manin connection,
on $\M$,\\
\item[(iv)] $\p$ is an odd homological (i.e.\ $\p^2=0$) vector field
on $\M$ such that $[E,\p]=\p$, $\nabla_X\nabla_Y\nabla_Z \p=0$ for any 
horizontal vector fields, $X,Y$ and $Z$ on $\M$, and  
$\p I\subset I^2$, $I$ being the ideal of the 
distinguished point in $\M$, \\
\item[(v)] $[\mu_n:\otimes^n \T \rar \T]$, $n\in {\Bbb N}$, is a 
homotopy class of smooth unital strong homotopy
commutative ($C_{\infty}$) algebras  defined on the 
tangent sheaf, $\T$, to $\M$, such that
$Lie_E \mu_n= \frac{1}{2} n\mu_n$, for all 
$n\in {\Bbb N}$, and  $\mu_1$ is given by
$$
\Ba{rccc}
\mu_1: & \T & \lon & \T \\
& X &\lon & \mu_1(X):=[\p,X].
\Ea
$$
\\
\item[(v)] $e$ is the flat unit, i.e.\ 
an even vector field on $\M$ such that $[\p, e]=0$, $\nabla e=0$,
$\mu_2(e,X)=X$, $\forall X\in \T$, and $\mu_n(\ldots, e,\ldots )=0$
for all $n\geq 3$.\\
\end{itemize}

\bip

\paragraph{\bf 3.6.1. Theorem. } {\em There is a canonical functor
from the category of pairs $(\fg, F)$, where $\fg$ is
$L_{\infty}$-formal unital homotopy Gerstenhaber algebra and $F$ a formality
map, to the category of  pre-${\mathit Frobenius}_{\infty}$ 
manifolds.}

\bip

\Proof The desired statement follows immediately from Theorem~2.7.1
and a version of Theorem-Construction~3.4.2 where the formality map
$F$ is used to transfer the $G_{\infty}$-structure from the algebra 
to its cohomology.
\hfill $\Box$

\bip

\paragraph{\bf 3.6.2. Theorem. } {\em If a homotopy Gerstenhaber algebra
$\fg$ is quasi-isomorphic, as a $L_{\infty}$-algebra, to an Abelian
dLie algebra, then the tangent sheaf, $\T$, to its cohomology
viewed as a linear supermanifold is canonically a sheaf
of unital graded commutative associative algebras.}

\bip

\Proof In this case $\p=0$ and $\mu_2$, which is now defined uniquely,
 makes $\T$ into a sheaf of  unital graded commutative 
associative algebras.
\hfill $\Box$

\bip


\section{Perturbative construction of $F_{\infty}$-invariants}


 The purpose of this section is to give second 
 ``down-to-earth'' proofs of some of the main claims
of this paper. Our approach here is 
based on perturbative solutions
of  algebro-differential equations rather than on the
 homotopy technique used in the two previous Sections.

\sip

First comes a perturbative proof of the Smoothness Theorem~ 2.5.6.
 
\bip

\sip

\paragraph{\bf 4.1. Theorem (Chen's construction).} 
{\em For any differential Lie superalgebra $\fg$, 
there exists a versal
element, $\Gamma\in k[[t]]\ot \fg$, and an odd derivation, 
$\p: k[[t]] \lon k[[t]]$,  such that $\p^2=0$ and the equation, 
$$
d\Gamma + \vec{\p}\Gamma + \frac{1}{2} [\Gamma\bullet \Gamma] =0
$$
holds.
Moreover, for any quasi-isomorphism  of complexes of vector spaces, 
$\phi: (\fg, d) \lon (\bH, 0)$, $\Gamma$  
 may be normalized so that $\phi(\Gamma_{[n]})=0$ for all $n\geq 2$.
}

\bip

We shall prove this Theorem by induction using (twice) the following
Lemma which is merely a trancated version of Remark~2.5.1.

\bip

\paragraph{\bf 4.1.1. Lemma.} {\em Assume the elements
$\Gamma_{(n)}=\sum_{k=0}^n \Gamma_{[k]}\in k[[t]]\ot 
\fg$ and  $\p_{(n)}=\sum_{k=0}^{n} \p_{[n]}\in {\mathsf Der}\,  k[[t]]$
satisfy
$$
d\Gamma_{(n)}+ \vec{\p}_{(n)} \Gamma_{(n)} + \frac{1}{2}\left[\Gamma_{(n)}
\bullet
\Gamma_{(n)}\right] = 0  \ \bmod I^{n+1}.
$$
Then 
$$
\psi_{[n+1]}:= d\Gamma_{(n)}+ \vec{\p}_{(n)} \Gamma_{(n)} + 
\frac{1}{2}\left[\Gamma_{(n)}\bullet\Gamma_{(n)}\right] \, \bmod I^{n+2}
$$
satisfies
$$
d \psi_{[n+1]} = - \vec{\p}_{(n)}^{\, 2} \left(\Gamma_{(n)}\right) \, 
\bmod I^{n+2}.
$$
}

\sip

\bip

 {\bf 4.1.2. Proof of the Theorem.}  
Let
$$
\phi: (\fg, d) \lon (\bH, 0),
$$
be a quasi-isomorphism,
i.e.\ a morphism of complexes 
inducing an isomorphism on cohomology.
Since $\fg$ is defined over a field, such a 
quasi-isomorphism always exists (note that we do not ask for any sort of a 
relationship between $\phi$ and the Lie brackets).

\sip

Let $e_i$ be any representatives of the cohomology classes $[e_i]$ in 
$\Ker\, d\subset \fg$. We may assume without loss of generality
that $\phi(e_i)=[e_i]$. Then choosing $\Gamma_{[0]}=0$, 
$\Gamma_{[1]}:=\sum_{i=1}^{p+q} t^i e_i$, and $\p_{[0]}=\p_{[1]}=0$
we get the data $(\Gamma_{(1)}, \p_{(1)})$ satisfying the Master equation 
modulo terms in $I^2$ and
the nilpotency condition $\p^2=0$ modulo terms in $I^3$.

\sip

Assume we have
constructed a versal element 
$\Gamma_{(n)}=\sum_{k=1}^n \Gamma_{[k]}\in k[[t]]\ot 
\fg$ and an odd vector field $\p_{(n)}=\sum_{k=2}^{n} \p_{[n]}$ on $\M$
such that the equations
$$
P_n: \left\{ \Ba {l}
d\Gamma_{(n)}+ \vec{\p}_{(n)} \Gamma_{(n)} + \frac{1}{2}\left[\Gamma_{(n)}
\bullet
\Gamma_{(n)}\right] = 0 \, \bmod I^{n+1}\vspace{3mm} \\
\p_{(n)}^2 = 0 \, \bmod I^{n+2}.\Ea \right.
$$
are satisfied.

\sip

Let us show that there exists $\Gamma_{[n+1]}\in k[[t]]\ot \fg$
and $\p_{[n+1]}\in H^0(\cT M_{\bH})$ such that 
$$
\Gamma_{(n+1)} = \Gamma_{(n)} + \Gamma_{[n+1]}, \ \ \ \ \ \ \ \
\p_{(n+1)} = \p_{(n)} + \p_{[n+1]}, 
$$
satisfy the equations $P_{n+1}$.

\sip

Note that, in the notations of  Lemma~4.4.1, one has 
$$
d\Gamma_{(n+1)}+ \vec{\p}_{(n+1)} 
\Gamma_{(n+1)} + \frac{1}{2}\left[\Gamma_{(n+1)}\bullet
\Gamma_{(n+1)}\right] \bmod I^{n+2} = 
d\Gamma_{[n+1]} + \psi_{[n+1]} + \vec{\p}_{[n+1]}\Gamma_{[1]}.
$$
Let us now  define $\vec{\p}_{[n+1]}$ by setting
$$
\vec{\p}_{[n+1]}\Gamma_{[1]}:= -\phi(\psi_{[n+1]}).
$$
As $d \psi_{[n+1]}=0$ by Lemma~4.4.1 and the second equation of $P_n$, 
we conclude that
$$
\psi_{[n+1]} + \vec{\p}_{[n+1]}\Gamma_{[1]} 
\in (\Ker\phi \cap \Ker d)\ot k[[t]]_{[n+1]}.
$$
Since $\phi$ is a quasi-isomorphism, $\Ker\,\phi \cap \Ker\, d =\Img\, d$.
Hence, there exists $\Gamma_{[n+1]}\in k[[t]]\ot\fg$ such that
$$
d\Gamma_{[n+1]} = - \psi_{[n+1]} - \vec{\p}_{[n+1]}\Gamma_{[1]}. 
$$
Thus the first equation of the system $P_{n+1}$ holds. This implies, 
by Lemma~2.4.2,
\Beqrn
d\psi_{[n+2]} &=& -\vec{\p}_{(n+1)}^{\, 2} \Gamma_{(n+1)} \, 
\bmod I^{n+3} \\
&=& -\vec{\p}_{(n+1)}^{\, 2} \Gamma_{[1]} \, \bmod I^{n+3}.
\Eeqrn
Applying $\phi$ to both sides of this equation, we get
$$
\vec{\p}_{(n+1)}^{\, 2} \phi(\Gamma_{[1]})=0
$$ 
implying the second equation of the system $P_{n+1}$,
$$
\vec{\p}_{(n+1)}^{\, 2}=0\, \bmod I^{n+3},
$$
and completing thus the inductive procedure.

\sip

Finally, we note that $\Ker d\, +\, \Ker \phi=\fg$ for $\phi$ is a 
quasi-isomorphism.
Hence we can always adjust $\Gamma_{[n+1]}$, $n\geq 1$,  so that it lies in
$\Ker\phi$. \hfill $\Box$.

\bip

\paragraph{\bf 4.1.3. Remarks.} {(i)} The role of $\p$ in the 
Chen's construction is to absorb
all the obstructions  so that constructing a versal solution to the 
Master equation poses no problem (cf.\ Smoothness Theorem 2.5.6).

\bip

(ii)\, The Chen's differential $\p$ is completely determined by $\Gamma$. 
Indeed, the Master equations imply,
$$
\vec{\p}\phi(\Gamma) = -\frac{1}{2}\phi\left([\Gamma\bullet\Gamma]\right).
$$
Decomposing,
$$
\phi(\Gamma)= \sum_i f^i(t) [e_i],
$$
we note that  $f^i(t)=t^i \bmod I^2$. Hence the functions $f^i(t)$ define a 
coordinate system
on $\M$ and the values,  ${\p} f^i(t)$, completely determine 
the differential $\p$.

\sip

In particular, if $\Gamma$ is $\phi$-{\em normalized}, i.e.\  
$\phi(\Gamma_{[n\geq 2]})=0$, then $\p$ can be computed by the formula
$$
\vec{\p}\left(\sum_{i=1}^{p+q} t^i [e_i]\right) = -\frac{1}{2}
\phi\left( [\Gamma \bullet \Gamma ]
\right).
$$

\bip

{(iii)}\, We shall understand from now on a versal solution, $\Gamma$, of 
the Master equations and the associated Chen differential $\p$ as, 
respectively, global sections of the sheaves $\fg\ot \f_{\M}$ and $\T$ 
on $\M$ (in practical terms,
this essentially fixes their transformation properties under  arbitrary 
changes of coordinates on the cohomology supermanifold). 

\sip

We  call sometimes $\Gamma$ 
 a {\em Master function}.

\bip

(iv) The  argument in (ii) also downplays the role of the 
quasi-isomorphism $\phi$ used in the Chen construction.
If  $\Gamma$ is normalised with respect to a quasi-isomorphism
$\phi:(\fg,d)\rar (\bH,0)$, then,  for any other quasi-isomorphism 
${\phi}'$, the {\em same}\, $\Gamma$  can be viewed as ${\phi}'$-normalized,  
but in a new coordinate
system ${t'}^i=f^i(t^j)$ given by ${\phi'}(\Gamma)= f^i(t^j)[e_i]$.
{ Thus varying  quasi-isomorphism $\phi$ used in the 
construction of $\Gamma$
amounts to varying flat structure on the pointed supermanifold $\M$}.

\bip

(v) Chen has actually invented his differential $\p$  in the context of 
differential associative algebras \cite{Chen}. Its Lie algebra analogue, 
Theorem~4.1, is  due to
Hain \cite{Hain}.

\sip

\bip

\paragraph{\bf 4.2. Gauge equivalence.} Let us consider
the following action, called a {\em gauge transformation}, 
of $\fg_{\tln}\ot I$ on $\fg\ot k[[t]]$:
$$
\begin{array}{ccc}
 \fg\ot I\, \times \, \fg\ot k[[t]] & \lon & \fg\ot k[[t]]\vspace{3mm} \\
 g \ot \Gamma & \lon &   \Gamma^g := e^{\ad_g}\Gamma -
\frac{e^{\ad_g}-1}{\ad_g}(d+\vec{\p})g.
\end{array}
$$

\bip

\paragraph{\bf 4.2.1. Lemma.} {\em If\, $\Gamma\in \fg\ot k[[t]]$ is a 
Master function,
then, for any $g\in \fg_{\tln}\ot I$, the function $\Gamma^g$
is also a Master function, and both these share the same Chen differential.}

\bip

\Proof We have to show that the equation
$d\Gamma + \ovr{\p}\Gamma + \frac{1}{2}[\Gamma\bullet \Gamma]=0$
implies
$$
d\Gamma^g + \ovr{\p}\Gamma^g + \frac{1}{2}[
\Gamma^g\bullet \Gamma^g]=0. 
$$
This follows immediately from the well-known formulae \cite{GM},
$$
e^{\ad_g}d e^{-\ad_g} = d - \ad_{\frac{e^{\ad_g}-1}{\ad_g}dg}, \ \ \ \ \ \ \ \
e^{\ad_g}\p e^{-\ad_g} = \p - \ad_{\frac{e^{\ad_g}-1}{\ad_g}\vec{\p} g},
\ \ \ \ \ \ \ \  e^{\ad_g}\ad_{\Gamma}e^{-\ad_g}=\ad_{e^{\ad_g}\Gamma}, 
$$
and 
$$
e^{\ad_g}\left[(\ldots)\, \bullet\, (\ldots)\right] = 
\left[e^{\ad_d}(\ldots)\, \bullet\, e^{\ad_g}(\ldots)\right].
$$
\hfill $\Box$

\bip

\paragraph{\bf 4.2.2. Theorem.} {\em Let\,  $\fg$ be a 
differential Lie algebra. For any two Master functions on $\M$,
 $\Gamma$ and\, $\Gamma'$, there is a gauge function
$g\in \Gamma(\M,  \fg\ot I)$ and a 
diffeomorphism $f:(\M,0)\rar(\M,0)$
such that
$\Gamma'= f^*(\Gamma^g)$  
and\,  $\p=f_*(\p')$.}

\bip

{\bf A sketch of the proof.} 
Let us fix a quasi-isomorphism $\phi:(\fg,d)\rar (\bH,0)$ of 
complexes of Abelian groups, and a coordinate system on $\M$ in which
$\Gamma'$ is $\phi$-normalized. 


\sip

We have $\Gamma'_{[1]}=\Gamma_{[1]}-dg_{[1]}$, for some 
$g_{[1]}\in \Gamma(\M,  \fg\ot I)$, and 
 $\p'_{[1]}=\p_{[1]}=0$. Hence
$\Gamma'=  \Gamma^{g_{[1]}} \ \bmod I^2$ and there is a unique
diffeomorphism $f_1:\M\rar\M$ such that the Master function 
$\Gamma'':=f_1^*(\Gamma^{g_{[1]}})$
is $\phi$-normalized and $\Gamma''_{[1]} = \Gamma'_{[1]}$. Hence,
$$
d(\Gamma'_{[2]} - \Gamma''_{[2]}) + 
\ovr{(\p'_{[2]} - \p''_{[2]})}\Gamma_{[1]}=0
$$
implying $\p'_{[2]}= \p''_{[2]}$ and  $d(\Gamma'_{[2]} - \Gamma''_{[2]})=0$.
Since $\phi(\Gamma'_{[2]} - \Gamma''_{[2]})=0$ and 
$\phi$ is a quasi-isomorphism,
there exists $g_{[2]}\in  \fg \ot I^2$ such that 
$\Gamma'_{[2]} - \Gamma_{[2]} = d g_{[2]}$. Hence
$$
\Gamma'=  (\Gamma'')^{g_{[2]}}= f_1^*(\Gamma^{g_{(2)}}) \ \bmod I^3.
$$

\sip

Continuing  by induction and using Lemma~4.2.1 
one easily obtains the desired result.
 \hfill $\Box$

\bip

\paragraph{\bf 4.2.3. Corollary.} {\em The Chen's vector field $\p$ on $\M$
is an invariant of $\fg$.}



\sip

\bip

\paragraph{\bf 4.3. Differential on $\T$.}
We fix from now on a dG-algebra $\fg$ and a Master function $\Gamma$ on $\M$.
The latter is not defined canonically, though the associated Chen differential
$\p$ {\em is}. 

\sip

We also  fix a  quasi-isomorphism,
$\phi: (\fg,d)\rar (\bH,0)$, of complexes of Abelian groups. This puts
 {\em no}\, restriction 
whatsoever on the dG-algebra under consideration. Moreover,
our main results will not depend
on the particular choices of $\Gamma$ and $\phi$ we have made  --- these two are 
no more than working tools.

\sip

The global vector field 
$\p$ on $\M$ makes $\T$ into a sheaf of complexes with the differential
$$
\begin{array}{rccc}
\dv: & \T & \lon & \T \\
& \X\,  & \lon & \dv \X:=[\p,\X\,],
\end{array}
$$
where $[\ ,\ ]$ stands for the usual commutator of (germs) of vector fields.
Indeed,
$$
(\dv)^2 \X\,=[\p,[\p,\X\,]]= \frac{1}{2}[[\p,\p],\X\,]=[\p^2,\X\,]=0,
$$
where we have used the Jacobi identity and the fact that $\p^2=0$.

\sip

This, of course, induces a differential on the sheaf of tensor products,
$\T^{\ot m} \ot (\T^*)^{\ot n}$ (and on the 
associated vector space of 
global sections),
which we denote by the same symbol $\dv$.

\bip

\paragraph{\bf 4.4. Deformed dG-algebra.} It is easy to check that  the map
$$
\begin{array}{rccc}
\dg: & k[[t]]\ot \fg & \lon & k[[t]]\ot \fg \vspace{3mm}\\
& a  & \lon & \dg a :=d a + \vec{\p}a  + 
[\Gamma\bullet a].
\end{array}
$$
satisfies
\Bi
\item[(i)] $(\dg)^2=0$,\\
\item[(ii)] $\dg(a\cdot b) =  
(\dg a)\cdot b + (-1)^{\tl{a}}a\cdot \dg b$ \\
\item[(iii)] $\dg\left[a\bullet b \right] =  \left[ \dg a\bullet b\right]
-(-1)^{\tl{a}}\left[a\bullet \dg b\right]$
\Ei
implying that the data
$
\left(k[[t]]\ot\fg, [\, \bullet\, ],\, \cdot \,, \dg\right)
$ 
is a dG-algebra.

\sip

The differentials $\dg$ and $\dv$ make the sheaf $\fg\ot (\T^*)^{\ot k}$
on $\M$ into a sheaf of complexes with the differential which 
we denote by $\Dg$.
For example, for any germ $\Phi\in \fg\ot \T^*$ and any germ $\X\in \T$ 
over the same point in $\M$,
$$
(\Dg \Phi)(\X\,) := \dg\Phi(\X\,) - (-1)^{\tl{\Phi}} \Phi(\dv \X\,).
$$
The vector space $\Hom(\T^{\ot k}, \fg\ot \f_{\M})$ is also a complex 
with the differential denoted by the same symbol $\Dg$.

\bip

\paragraph{\bf 4.5. Morphism of sheaves of complexes.}
The versal solution $\Gamma$ gives rise to a morphism of $\f_{\M}$-modules,
$$
\begin{array}{rccc}
\Gai: & \cT_{\bH} & \lon & \f_{\M}\ot \fg \vspace{3mm}\\
& \X\,  & \lon & \Gai(\X\,):= \ovr{\X\,\,}\Gamma.
\end{array}
$$
It is not hard to check that $\Gai$ is a monomorphism.

\bip
\paragraph{\bf 4.5.1. Lemma.} {\em The element}\, $\Gai \in \Hom(\T, \fg\ot \f_{\M})$
{\em is cyclic, i.e.}
$$
\Dg \Gai = 0.
$$

\Proof  Applying ${\X\,}\in  \T$ to both sides of the equation
$$
d\Gamma + \vec{\p}\Gamma + \frac{1}{2} [\Gamma\bullet \Gamma] =0
$$
we get
$$
(-1)^{\tl{X}} \dg(\ovr{\X\,\,}\Gamma) + \ovr{[\X,\p]\,\,}\,\Gamma =0,
$$
implying $(\Dg\Gai)(\X\,)=0$. \hfill $\Box$

\bip

\paragraph{\bf 4.5.2. Corollary.} 
{\em For any}\, $\X\,\in \T$, \ 
$
\dg(\ovr{\X\,\,}\Gamma) = \ovr{\dv \X\,\,}\,\Gamma$.

\bip

\paragraph{\bf 4.5.3 Corollary.} {\em For any}\, $\chi
\in \Hom(\ot^k \T, \T)$ {\em
one has}
$$
\Dg(\Gai\circ \chi) = \Gai\circ (\dv \chi).   
$$
\sip

\Proof We have, using Corollary~4.5.2,
\Beqrn
\Dg(\Gai\circ \chi)(\X_1, \ldots, \X_k)&=&
\dg \left(\ovr{\chi(\X_1, \ldots, \X_k)}\Gamma\right) -
(-1)^{\tl{\chi}}
\ovr{\chi(\dv \X_1, \ldots, \X_k)}\Gamma\\ && -\, \dots \,
- (-1)^{\tl{\chi}+ \tl{X}_1+\ldots + \tl{X}_{k-1}}
\ovr{\chi(\X_1, \ldots, \dv \X_k)}\Gamma\\
&=& \left(\ovr{\dv\chi(\X_1, \ldots, \X_k)}\Gamma\right) -
(-1)^{\tl{\chi}}
\ovr{\chi(\dv X_1, \ldots, X_k)}\Gamma\\ && -\, \dots \,
- (-1)^{\tl{\chi}+ \tl{X}_1+\ldots + \tl{X}_{k-1}}
\ovr{\chi(\X_1, \ldots, \dv \X_k)}\Gamma\\
&=& \left(\Gai\circ \dv \chi\right) (\X_1, \ldots, \X_k)  
\Eeqrn
for arbitrary $\X_1,\ldots, \X_k\in \T$. \hfill $\Box$


\sip

Therefore, the map
$$
\Gai: \left(\T, \dv\right) \rar \left(\fg\ot \f_{\M}, \dg\right)
$$
is a morphism of sheaves of complexes. Note that the ``projection'' map $s\circ \phi:
\fg\ot \f_{\M}\rar \T$ satisfies
$s\circ\phi\circ\Gai = \Id$\, but does {\em not}, in general, 
respect the differentials.

\sip

Analogously one shows that the morphism $\Gai\cdot \Gai$ defined by the commutative
diagram 
$$
\diagram
\T\ot\T \drto_{\Gai\cdot\Gai}  \rto^{\Gai\ot \Gai} & \fg\ot\fg\ot\f_{\M}
\dto^{\,\cdot\ot Id} \\
& \fg\ot\f_{\M}
\enddiagram
$$
defines a cyclic element in $(\Hom(\T^{\ot 2},\fg\ot \f_{\M}), \Dg)$. In a 
similar way one 
uses muliplicative structure in $\fg$ to construct cyclic elements
$\Gai\cdot\Gai\cdot\Gai$ etc.\footnote{The cyclicity of 
$\Gai\cdot \Gai$, etc., relies
 on the Poisson identity holding in 
$(\fg,[\, \bullet \, ], \, \cdot\,, d )$.}

\sip

For future reference we define a morphism $\Gai_{(n)}\in \Hom(\cT,\fg\ot\f_{\M})$
by setting  $\Gai_{(n+1)}(\X\,):= \ovr{\X\,} \Gamma_{(n+1)}$. Similarly one defines
$\Gai_{(n)}\cdot \Gai_{(n)}$, etc.

\bip

\paragraph{\bf 4.6. Multiplicative structure in $\T$.} We will show 
in this subsection that, for any dG-algebra $\fg$,
the associated tangent sheaf $\T$ is always a sheaf of 
differential associative algebras (defined uniquely up to a homotopy).

\bip

\paragraph{\bf 4.6.1. Theorem.} {\em
There exists an even  morphism of sheaves}, $\mu\in \Hom(\T^{\ot 2}, \T)$,
{\em such that $\dv\mu=0$  and
 the diagram
$$
\diagram
\T\ot\T \drto_{\Gai\cdot\Gai}  \rto^{\mu} & \T
\dto^{\Gai} \\
& \fg\ot\f_{\M}
\enddiagram
$$
is commutative at the cohomology level, i.e.
$$
[\Gai\cdot \Gai] = [\Gai\circ \mu]
$$
in the cohomology sheaf}\, $\Ker \Dg/ \Img \Dg$
{\em associated with the sheaf of complexes}\, 
$( Hom(\T^{\ot 2},\fg\ot \f_{\M}), \Dg)$.

\bip

\Proof We have to show that there exists  $\mu\in \Hom(\T^{\ot 2}, \T)$
such that
\Beq
\dv \mu(\X\,,\Y\,) = \mu(\dv \X\, , \Y\,) +
(-1)^{\tl{X}}\mu(\X\,, \dv \Y\,)
\label{1}
\Eeq
and
\Beq
\ovr{\X\,\,}\Gamma\cdot \ovr{\Y\,\,}\Gamma = 
\ovr{\mu(\X\, ,\Y\,)\,}\, \Gamma + (\Dg A)(\X\,,\Y\,)
\label{2}
\Eeq
for some $A\in \Hom(\T^{\ot 2},\fg\ot \f_{\M})$ and 
 any $\X\,,\Y\,\in \T$. We shall proceed 
by induction and assume, without loss of generality,
that the vector fields $\X\,$ and $\Y\,$ are constant, i.e.\ $\nabla \X=\nabla \Y=0$.

\sip

The above equations 
can obviously be satisfied $\bmod I$: just set
$$
\mu_{[0]}(\X\,,\Y\,):= \phi(\ovr{\X\,\,}\Gamma_{[1]}\cdot \ovr{\Y\,\,}\Gamma_{[1]}).
$$ 
Indeed,  
$$
\ovr{\X\,\,}\Gamma_{[1]}\cdot \ovr{\Y\,\,}\Gamma_{[1]} - 
\ovr{\mu_{ [0]}(\X\,,\Y\,\,)}\, \Gamma_{[1]} \in \Ker\, \phi \cap \Ker\, d
$$
and hence this expression is $d$-exact. Denote it by $d A_{[0]}(\X\,,\Y\,)$. 
(We can always normalise
$A_{[0]}$ so that it lies in $\Hom(\T^{\ot 2}, \ker\phi\ot \f_{\M})$.)
This solves (\ref{2}) $\bmod I$. The equation (\ref{1})$\bmod I$  is trivial
(recall that $\p_{[<2]}=0$).

\sip

Assume now that we have constructed $\mu_{(n)}$ and $A_{(n)}$ so that
the equations 
\Beq
\dv_{(n)} \mu_{(n-1)}(\X\,,\Y\,) = \mu_{ (n-1)}({\dv}_{(n)} \X\,, \Y\,) +
(-1)^{\tl{X\,}}\mu_{ (n-1)}(\X\,, {\dv}_{(n)} \Y\,) \ \bmod I^{n+1}
\label{3}
\Eeq
\Beq
\ovr{\X\,\,}\Gamma_{(n+1)}\cdot \ovr{\Y\,\,}\Gamma_{(n+1)} = 
\ovr{\mu_{ (n)}(\X\,,\Y\,)\,}\, \Gamma_{(n+1)} + ({\Dg}_{(n)} A_{(n)})(\X\,,\Y\,)
\ \bmod I^{n+1}
\label{4}
\Eeq
hold.

\sip

The Theorem will be proved if we  find $\mu_{[n+1]}$ and $A_{[n+1]}$ satisfying
$$
\dv_{(n+1)} \left(\mu_{(n)}(\X\,,\Y\,) + \mu_{[n+1])}(\X\,,\Y\,)\right)
= \mu_{ (n)}({\dv}_{(n+1)} \X\,, \Y\,) +
(-1)^{\tl{X}}\mu_{ (n)}(\X\,, {\dv}_{(n+1)} \Y\,) \ \bmod I^{n+2}
$$
and
$$
\ovr{\X\,\,}\Gamma_{(n+2)}\cdot \ovr{\Y\,\,}\Gamma_{(n+2)} - 
\ovr{\mu_{ (n)}(\X\,,\Y\,)\,}\, \Gamma_{(n+2)} - 
\ovr{\mu_{ [n+1]}(\X\,,\Y\,)\,}\, \Gamma_{[1]} 
 - ({\Dg}_{(n+1)} A_{(n)})(\X\,,\Y\,)
= \ \ \ \ \ \ \ \ \ \ \ \
$$
$$
\ \ \ \ \ \ \ \ \ \ \ \ \ \ \ \ \ \ \ \ \ \ \ \ \ \ \ \  \ \ \ \ \ \ \ \ =
d A_{[n+1]}(\X\,,\Y\,)
\ \bmod I^{n+2}
$$

Defining
$$
\mu_{ [n+1]}(\X\,,\Y\,):= \phi\left(
\ovr{\X\,\,}\Gamma_{(n+2)}\cdot \ovr{\Y\,\,}\Gamma_{(n+2)} - 
\ovr{\mu_{ (n)}(\X\,,\Y\,)\,}\, \Gamma_{(n+2)} - ({\Dg}_{(n+1)} A_{(n)})(\X\,,\Y\,)\right)
$$
we ensure that the morphism
$$
\lambda_{[n+1]}(\X\,,\Y\,):= \left(
\Gai_{(n+1)}\cdot \Gai_{(n+1)} - \Gai_{(n+1)}\circ \mu_{(n)} -
{\Dg}_{(n+1)} A_{(n)}) - \Gai_{(0)}\circ \mu_{[n+1]}\right)(\X\,,\Y\,) \ \ \bmod I^{n+2}
$$
take values in the sheaf in $\Ker\phi\ot \f_{\M}$.
Since it vanishes modulo $I^{n+1}$, we have, modulo $I^{n+2}$,
\Beqrn
d \lambda_{[n+1]}(\X\,,\Y\,) &=& (\Dg_{(n+1)} \lambda_{[n+1]})(\X\,,\Y\,)\\
&=& - \left(\Dg_{(n+1)}(\Gai_{(n+1)}\circ \mu_{(n)})\right)(\X\,,\Y\,) \ \\
&=&  -\dg_{n+1}(\ovr{\mu_{(n)}(\X\,,\Y\,)\,}\Gamma) + 
\ovr{\mu_{(n)}(\dv_{(n+1)}\X\,, \Y\,)\,}\Gamma
 + (-1)^{\tl{X}}\, \ovr{\mu_{(n)}(\X\,,\dv_{(n+1)}\Y\,)\,}\Gamma  \\
&=& - \ovr{(\dv_{(n+1)} \mu_{(n)}(\X\,,\Y\,) - \mu_{(n)}(\dv_{(n+1)}\X\,, \Y\,) -
(-1)^{\tl{X}} \mu_{(n)}(\X\,,\dv_{(n+1)}\Y\,))}\Gamma 
\Eeqrn
where we have used Corollary~4.5.2 and the fact that $\Dg(\Gai\cdot \Gai)=0$.
Applying $\phi$ to the last equation, we get
$$
\dv_{(n+1)} \mu_{(n)}(\X\,,\Y\,) = \mu_{(n)}(\dv_{(n+1)}\X\,, \Y\,) +
(-1)^{\tl{X}} \mu_{(n)}(\X\,,\dv_{(n+1)}\Y\,) \ \  \bmod I^{n+2}
$$
and hence
$$
d \lambda_{[n+1]}(\X\,,\Y\,)=0.
$$
Since $\Ker \phi\cap \Ker d\subset \Img d$, 
there exists $A_{[n+1]}(\X\,,\Y\,)$
(which can be chosen to lie in $\Ker\phi\ot \f_{\M}$)
such that
$$
\lambda_{[n+1]}(\X\,,\Y\,)= d A_{[n+1]}(\X\,,\Y\,).
$$
This completes the inductive procedure and hence the proof of the Theorem.
 \hfill  $\Box$

\bip

\paragraph{\bf 4.6.2. Definition.} An even morphism of sheaves, 
$\mu\in \Hom(\ot^2 \T,\T)$,
satisfying the conditions of Theorem~4.6.1 is called {\em induced}. 
The associated 
data $(\T,\delta,\mu)$ is called a sheaf of 
{\em induced differential algebras}.

\sip

Clearly, an induced product on $\T$ is supercommutative
if the product $\cdot$ in $\fg$ is supercommutative. 

\bip

\paragraph{\bf 4.7. (Non)Uniqueness.} How unique is the product $\mu$
induced on 
the tangent sheaf $\T$ by Theorem~4.6.1?
When is it associative? To address these questions 
we shall need the following
technical result.

\bip


\paragraph{\bf 4.7.1. Lemma.} {\em If}\, $\tau
\in \Hom(\ot^k \T, \T)$ and
$B \in \Hom(\ot^k \T, \fg\ot \f_{\M})$ 
{\em
satisfy the equation
$$
\Gai\circ \tau = \Dg B 
$$
then there exists}\, $\chi\in \Hom(\ot^k \T, \T)$ and
$C \in \Hom(\ot^k \T, \Ker\phi \ot \f_{\M})$ {\em such that
\Bi
\item[{\em (1)}]
$
B= \Gai\circ \chi + \Dg C$,
\item[{\em (2)}] $\tau=\dv \chi$, i.e.
\Beqrn
\tau(\X_1, \ldots, \X_k)
&=& \dv \chi(\X_1, \ldots, \X_k) - (-1)^{\tl{\chi}}
\chi(\dv \X_1, \ldots, \X_k)\\ && -\, \dots \,
- (-1)^{\tl{\chi}+ \tl{X}_1+\ldots + \tl{X}_{k-1}}
\chi(\X_1, \ldots, \dv \X_k)
\Eeqrn
for any $\X_1,\ldots, \X_k\in \T$.
\Ei
}

\bip

\Proof  Without loss of generality we may assume that (germs of) vectors
fields $\X_1, \ldots, \X_k$ are constant. In view of Corollary~4.5.3 
and injectivity of the map $\Gai$, it 
is enough to show that the equation
\Beq
\ovr{\tau(\X_1, \ldots, \X_k)}\Gamma = (\Dg B)(\X_1, \ldots, \X_k)\label{5}
\Eeq
implies
$$
B(\X_1, \ldots, \X_k)= \ovr{\chi(\X_1, \ldots, \X_k)}\Gamma 
+ (\Dg C)(\X_1, \ldots, \X_k)
$$
for some  $\chi\in \Hom(\ot^k \T, \T)$ and
$C \in \Hom(\ot^k \T, \Ker\phi \ot \f_{\M})$. We shall proceed by induction.

\sip

The equation (\ref{5})~$\bmod I$ is
$$
\ovr{\tau_{[0]} (\X_1, \ldots, \X_k)}\Gamma_{[1]} = dB_{[0]}(\X_1, \ldots, \X_k).
$$
Hence $\tau_{[0]}=0$ and $dB_{[0]}=0$. Set 
$$
\chi_{[0]}(\X_1,\ldots, \X_k):= \phi\left(B_{[0]}(\X_1,\ldots,\X_k)\right).
$$
Then $
B_{[0]}(\X_1,\ldots,\X_k) - \ovr{\chi_{[0]}(\X_1,\ldots, \X_k)}\Gamma_{[1]}
$
lies in $(\Ker\,\phi \cap \Ker\, d)\ot \f_{\M}$ and hence equals 
$dC_{[0]}(\X_1,\ldots,\X_k)$ for some  
$C_{[0]} \in \Hom(\ot^k \T, \Ker\phi \ot \f_{\M})$.

\sip

Assume that $\chi_{(n)}$ and $C_{(n)}$ are constructed so that the equations
$$
B_{(n)}(\X_1, \ldots, \X_k)= \ovr{\chi_{(n)}(\X_1, \ldots, \X_k)}\Gamma_{(n+1)}
+ (\Dg_{(n)} C_{(n)})(\X_1, \ldots, \X_k)\ \bmod I^{n+1}
$$
holds.

\sip

Let us show that there exist $\chi_{[n+1]}$ and $C_{[n+1]}$ satisfying
\Beqrn
B_{(n+1)}(\X_1, \ldots, \X_k)&=& \ovr{\chi_{(n)}(\X_1, \ldots, 
\X_k)}\Gamma_{(n+2)}
+ \ovr{\chi_{[n+1]}(\X_1, \ldots, \X_k)}\Gamma_{[1]}\\
&&  +\, (\Dg_{(n+1)} C_{(n)})(\X_1, \ldots, \X_k)
+ dC_{[n+1]} (\X_1, \ldots, \X_k)
\ \bmod I^{n+2},
\Eeqrn
or, equivalently, satisfying
$$
dC_{[n+1]}(\X_1, \ldots, \X_k)= \psi_{[n+1]}(\X_1, \ldots, \X_k) -  
\ovr{\chi_{[n+1]}(\X_1, \ldots, \X_k)}\Gamma_{[1]} \ \bmod I^{n+2},
$$
where we have set
$$
\psi_{[n+1]}(\X_1, \ldots, \X_k) := B_{(n+1)}(\X_1, \ldots, \X_k) -
\ovr{\chi_{(n)}(\X_1, \ldots, \X_k)}\Gamma_{(n+2)} -
 (\Dg_{(n+1)} C_{(n)})(\X_1, \ldots, \X_k).
$$
Since $\psi_{[n+1]}(\X_1, \ldots, \X_k)$ vanishes $\bmod I^{n+1}$, 
this is a monom of degree
$n+1$ and $t^i$, and hence, modulo $I^{n+2}$,
\Beqrn
d \psi_{[n+1]}(\X_1, \ldots, \X_k) &=& (\Dg_{(n+1)}\psi_{[n+1]})(\X_1, 
\ldots, \X_k)
 \\
&=& (\Dg_{(n+1)} B_{(n+1)}) (\X_1, \ldots, \X_k)  - 
\Dg_{(n+1)}(\Gai_{(n+1)}\circ 
\chi_{(n)}) (\X_1, \ldots, \X_k)    \\
&=& \ovr{\left(\tau_{(n+1)} (\X_1, \ldots, \X_k) - (\dv_{(n+1)} \chi_{(n)}) 
(\X_1, \ldots, \X_k)
  \right)}\,\Ga .\\
\Eeqrn
Applying $\phi$ to both sides of these equations we get
$$
\tau_{(n+1)}(\X_1, \ldots, \X_k) = (\dv_{(n+1)} \chi_{(n)}) 
(\X_1, \ldots, \X_k)\ \bmod I^{n+1}
$$ 
and hence
$$
d \psi_{[n+1]}(\X_1, \ldots, \X_k)=0 \ \bmod I^{(n+1)},
$$

\sip

We define $\chi_{[n+1]}$ by
$$
\chi_{[n+1]}(\X_1, \ldots, \X_k) := \phi( \psi_{[n+1]}(\X_1, \ldots, \X_k)).
$$
Then $\psi_{[n+1]}(\X_1, \ldots, \X_k) -  
\ovr{\chi_{[n+1]}(\X_1, \ldots, \X_k)}\Gamma_{[1]}\in \Ker\, d\cap \Ker\, \phi
\subset \Img d$. This proves the existence of $C_{[n+1]}$ and hence completes
the proof of the theorem. \hfill $\Box$

\bip

\bip

If $\mu', \mu''\in \Hom(\ot^2\T, \T)$ are two products as in Theorem~4.6.1, 
then
$$
\Gai\circ (\mu'-\mu'')= \Dg B
$$
for some $B\in \Hom(\ot^2 \T, \fg\ot \f_{\M})$, and hence, by Lemma~4.7.1,
$$
\mu'-\mu''= \delta \chi
$$ 
for some odd $\chi\in \Hom(\ot^2\T, \T)$, i.e.\ $\mu'$ and $\mu''$ are 
what is called
{\em homotopy equivalent}. 

\sip

On the other hand, if $\mu''$ is a product with the properties 
stated by Theorem~4.6.1,
then, for any odd $\chi\in \Hom(\ot^2\T, \T)$, the product
$$
\mu':= \mu'' + \delta \chi
$$
also enjoyes the properties of Theorem~4.6.1. 
Indeed, by Corollary~4.6.3,
\Beqrn
\Gai\circ \mu'' &=& \Gai\circ\mu' + \Gai\circ (\delta\chi) \\
&=& \Gai\circ\mu' + \Dg(\Gai\circ \chi)  
\Eeqrn
and hence
$$
[\Gai\cdot \Gai] = [\Gai\circ \mu'] = [\Gai\circ \mu']
$$
in the cohomology sheaf  $\Ker \Dg/ \Img \Dg$.

\sip

Thus the set of products $\mu$ induced on $\T$ by Theorem~4.6.1 
is a  principal homogeneous space over the 
Abelian group $\delta\Hom_{\tln}(\ot^2\T, \T)$. Hence all induced 
products on  each stalk of $\T$  combine into
a single homotopy class which we call {\em induced}. 

\bip

\paragraph{\bf 4.7.2. Theorem.} {\em The sheaf $\T$ is canonically 
a sheaf of induced
homotopy  classes of  differential algebras}.

\bip

\Proof  We have to show that the homotopy class of products induced
on $\M$ is an invariant of the dG-algebra under consideration, i.e.\
that
it is independent of the choice of a quasi-isomorphism $\phi$
and on the choice of a  Master function $\Gamma$ used in 
its construction. In view of  Remark~4.1.3(iv), it is enough to check
the invariance of the product under the gauge transformations,
$$
\Gamma \lon 
\Gamma^g := e^{\ad_g}\Gamma -
\frac{e^{\ad_g}-1}{\ad_g}(d+\vec{\p})g, \ \ \ \
g\in \Gamma(\M, \fg_{\tln}\ot \M).
$$
A straightforward analysis of the basic equation (\ref{1}) shows
that gauge transformation changes the tensor $A$,
$$
A^g= e^{\ad_g}\left(A(X,Y) - G\cdot\Gai - \Gai\cdot G
+ G\cdot \Dg G + G\circ \mu\right),
$$
where $G\in \Hom(\T, \fg\ot \f_{\M})$ is given by
$$
G(X):= (-1)^{\tl{X}}\frac{e^{\ad_g}-1}{\ad_g}(d+\vec{\p})\ovr{\X\,}g, 
$$
but leaves the product invariant, $\mu^g=\mu$. \hfill $\Box$

\bip

\sip

\paragraph{\bf 4.8. Identity in $\fg$ $\Rightarrow$ identity in $\T$.} 
If the dG-algebra under consideration, $\fg$, 
has the identity $e_0$, and the 
versal solution $\Gamma$ is approprietly normalized 
(see Remark 3.1.1), then
\Beqrn
\ovr{\delta (e)}\, \Gamma&=& \ovr{[\p, e]}\, \Gamma \\
&=& \vec{\p} e_0 + \vec{e}(d\Gamma + \frac{1}{2}[\Gamma\bullet \Gamma])\\
&=& 0 + de_0 + [e_0\bullet \Gamma]\\ 
&=& 0,
\Eeqrn
so that $\delta(e)=0$.
We shall show next that the induced
homotopy class of differential algebras on each stalk of $\T$  
containes a canonical subclass of {\em unital}\, differential algebras.

\bip

\paragraph{\bf 4.8.1. Theorem.} {\em If $\fg$ has the identity $e_o$, then 
$\T$ is canonically a sheaf of induced  homotopy classes of differential  
algebras with the identity $e$.}

\bip

{\bf A sketch of the proof.} It is enough to show that there exists
a $\delta$-closed element, $\mu \in \Hom_{\tlo}(\ot^2 \T, \T)$, such that,
for arbitrary (constant) $\X\,,\Y\,\in \T$, 
the equation  
$$
\ovr{\X\,}\Gamma\cdot \ovr{\Y\,}\Gamma = \ovr{\mu(\X\,,\Y\,)}\, \Gamma + (\Dg A)(\X\,,\Y\,)
$$
holds
for some $A\in \Hom(\T^{\ot 2},\fg\ot \f_{\M})$ satisfying 
$A(\X\,,e)=\X\,$ and $A(e,\Y\,)=\Y\,$ (cf.\ (\ref{2})).

\sip

Recall (see the proof  of Theorem~4.6.1) that at the lowest order we have,
\Beqrn
\mu_{[0]}(\X\,,\Y\,) &=& \phi(\ovr{\X\,}\Gamma_{[1]}\cdot \ovr{\Y\,}\Gamma_{[1]})\\
dA_{[0]}(\X\,,\Y\,)&=&\ovr{\X\,}\Gamma_{[1]}\cdot \ovr{\Y\,}\Gamma_{[1]} - 
\ovr{\mu_{ [0]}(\X\,,\Y\,)}\, \Gamma_{[1]}.
\Eeqrn
and 
hence $\mu_{[0]}(\X\,,e)=\mu_{[0]}(e,\X\,)=\X\,$ and $dA_{[0]}(\X\,,e)=dA_{[0]}(e,\X\,)=0$. We claim
that $A_{[0]}$ can be chosen to satisfy $A_{[0]}(\X\,,e)=A_{[0]}(e,\X\,)=0$. This can be 
achieved by a replacement,
$$
\Ba{ccccl}
A_{[0]}(\X\,,\Y\,) & \rar & A_{[0]}'(\X\,,\Y\,) & := & A_{[0]}(\X\,,\Y\,) - 
A_{[0]}(\X\,,e) \cdot \ovr{\Y\,}\Gamma - (-1)^{\tl{X}}\ovr{\X\,}\Gamma \cdot A_{[0]}(e,\Y\,) \\
&&&& +\, 
(-1)^{\tl{X}}\ovr{\X\,}\Gamma \cdot A_{[0]}(e,e) \ovr{\Y\,}\Gamma ,
\Ea
$$
which satisfies,
$$
dA'_{[0]}(\X\,,\Y\,)=dA_{[0]}(\X\,,\Y\,), \ \ \ \ A'_{[0]}(\X\,,e)= A'_{[0]}(e,\X\,)=0.
$$

\sip

This observation allows us to include into the inductive procedure of the proof
of Theorem~4.6.1 the additional assumptions
$$
\mu_{(n)}(\X\,,e)=\mu_{(n)}(e,\X\,)=\X\,, \ \ \ \ \ A_{(n)}(\X\,,e)= A_{(n)}(e,\X\,)=0,
$$
and show, by exactly the same argument as in the case $n=0$ above,
that they hold true for $n+1$. 

\sip

Thus there does exist a product $\mu$ from the induced homotopy class
satisfying $\mu(\X\,,e)=\mu(e,\X\,)=\X\,$. It is defined uniquely up to a transformation
$$
\mu \lon \mu +\delta \chi
$$
with $\chi$ satisfying $\chi(\X\,,e)=\chi(e,\X\,)=0$ for 
arbitrary $\X\,\in \T$.
Thus what is well-defined is the induced homotopy class of {\em unital}\, 
differential algebras. 
\hfill $\Box$

\bip

\sip


\paragraph{\bf 4.9. Theorem.} {\em For any 
(unital) dG-algebra $\fg$,  
the tangent sheaf $\T$ to its cohomology supermanifold
is canonically a sheaf of homotopy classes of (unital) 
$A_{\infty}$-algebras with 
\begin{itemize}
\item[(i)] $\mu_1= [\p, \ldots]$, $\p$ being the Chen's vector field, and \\
\item[(ii)] the homotopy class of $\mu_2$ being the induced homotopy class 
as in Theorem~4.6.1.
\end{itemize}
}

\bip

{\bf A sketch of the proof.} By Theorem~4.6.1, there 
exists a  product 
$\mu_2\in \Hom(\ot^2\T,\T)$ satisfying the equation
$$
\ovr{\mu_2(\X_1,\X_2)}\Gamma = \ovr{\X_1}\Gamma\cdot \ovr{\X_2}\Gamma +
(\Dg A_2)(\X_1, \X_2)
$$
for some odd $A_2\in \Hom(\ot^2 \T, \fg\ot \f_{\M})$ and arbitrary 
$\X_1,\X_2\in \T$.
We have, in the notations of subsection~3.2.3,
\Beqrn
\ovr{\Lambda_3(\X_1,\X_2,\X_3)}\Gamma &=& 
\ovr{\mu_2(\X_1,\mu_2(\X_2,\X_3)) - \mu_2(\mu_2(\X_1,\X_2),\X_3)}\Gamma \\
&=&  \ovr{\X_1}\Gamma\cdot \left(\ovr{\X_2}\Gamma \cdot \ovr{\X_3}\Gamma\right)
+ (\Dg A_2)(\X_1, \mu_2(\X_2,\X_4))  + \Gai(\X_1)\cdot (\Dg A_2)(\X_2,\X_3) 
\\ 
&&\hspace{-3mm} - 
\left(\ovr{\X_1}\Gamma\cdot \ovr{\X_2}\Gamma\right) \cdot \ovr{\X_3}\Gamma 
 - (\Dg A_2)(\mu_2(\X_1,\X_2),\X_3) - (\Dg A_2)(\X_1,\X_2)\cdot \Gai(\X_3) \\
&=& (\Dg B_3)(\X_1,\X_2,\X_3), 
\Eeqrn
where
\Beqrn
B_3(\X_1,\X_2,\X_3) &:=& (-1)^{\tX_1} \Gai(\X_1)\cdot A_2(\X_2,\X_3) - 
A_2(\X_1,\X_2)\cdot \Gai(\X_3) \\ 
&& +\, A_2(\X_1,\mu_2(\X_2,\X_3)) 
 - A_2(\mu_2(\X_1,\X_2),\X_3).
\Eeqrn
Here we used associativity of the dot product in $\fg$, $\Dg$-closedness of
$\Gai$ and $\delta$-closedness of $\mu_2$.

\sip

By Lemma~4.7.1, there exists $\mu_3\in \Hom_{\tlo}(\ot^3\T,\T)$ such that
the 3rd order associativity condition, $\Lambda_3 = \delta \mu_3$,
is satisfied, and
$$
\ovr{\mu_3(\X_1,\X_2,\X_3)}\Gamma= B_3(\X_1,\X_2,\X_3) + (\Dg A_3)(\X_1,\X_2,\X_3)
$$
for some $A_3\in \Hom_{\tlo}(\ot^3\T,\fg\ot \f_{\M})$.

\sip

Exactly the same procedure constructs inductively all the higher order products 
$\mu_n\in \Hom_{\tl{n}}(\ot^n \T,\T)$ 
satisfying the higher  order associativity conditions:
\begin{itemize}
\item[Step 1.] 
Assume that we have  constructed 
$\mu_k\in \Hom_{\tl{k}}(\ot^k \T,\T)$ and $A_k\in  
\Hom_{\tl{k}+\tln}(\ot^k\T,\fg\ot \f_{\M})$, $k=2,\ldots, n-1$, such that 
$\Lambda_k=\delta\mu_k$ ($k$-th order associativity condition) and
\Beqrn
\ovr{\mu_k(\X_1,\ldots,\X_k)}\Gamma&=& \sum_{i+j=k}(-1)^{(j+1)(\tX_1+\ldots+\tX_i)+i+1}
A_i(\X_1,\ldots,\X_i)\cdot A_i(\X_{i+1}, \ldots, \X_{k}) \\
&&\hspace{-7mm}  + \hspace{-1mm} \sum_{{i+j=k+1\atop i\geq 2}\atop j\geq 2}\sum_{l=0}^{i-1} 
(-1)^r A_i(\X_1, \ldots, \X_l, 
\mu_j(\X_{l+1}, \ldots, \X_{l+j}), \X_{l+j+1}, \ldots, \X_{k})\\
&&\hspace{-7mm} + \ (\Dg A_k)(\X_1,\ldots,\X_k) \vspace{3mm}\\
&=:& B_k(\X_1,\ldots, X_k) +  (\Dg A_k)(\X_1,\ldots,\X_k), 
\Eeqrn
where we have set $A_1:=\Gai$ and $r=\tl{j}(\tl{\X_1}+\ldots \tl{\X_1}) + 
\tl{l}(\tl{j}-1)
+ (\tl{i}-\tln)\tl{j} +1$. \\

\item[Step 2.] Use the above expressions for 
$\Gai\circ \mu_k$, $k=2,\ldots,n-1$,
to show that
$$
\ovr{\Lambda_n(\X_1,\ldots,\X_n)}\Gamma = (\Dg B_n) (\X_1,\ldots,\X_n).
$$

\item[Step 3.] Apply Lemma~3.6.1 to conclude that there exists 
$\mu_n$ such that 
$\Lambda_n=\delta\mu_n$
($n$-th order associativity condition) and
$\Gai\circ \mu_n= B_n + \Dg A_n$ for some 
$A_n\in \Hom_{\tl{n}+\tln}(\ot^k\T,\fg\ot \f_{\M})$

\end{itemize}

\sip

Finally, we note that at each stage of the above construction the $n$th product
$\mu_n$ is defined only up to a $\delta$-exact term, $\delta f_n$. 
These arbitrary  terms combine all together into a homotopy of the 
$A_{\infty}$-structure
$(\T, \mu_*)$. \hfill $\Box$

\bip

\paragraph{\bf 4.9.1. Corollary.} {\em The cohomology sheaf on}\,  $\M$,
$$
\Ho := \frac{\Ker\, \delta}{\Img\, \delta},
$$
{\em is canonically  a sheaf of induced (unital) associative algebras}.

\bip

\paragraph{\bf 4.9.2. Corollary.} {\em 
The tangent sheaf, $\cT\cM_{\mathrm smooth}$,
 to the smooth part of the extended
Kuranishi space is canonically a sheaf of
induced (unital) associative algebras}.

\bip

\paragraph{\bf 4.9.3. The Euler field.} If the dG-algebra $\fg$ 
under consideration is $\Z$-graded, then the  cohomology $\bH$ and hence 
its dual $\bH^*$ are also $\Z$-graded. We make $k[[t]]$ into a $\Z$-graded 
ring my setting
$$
k[[t]]=\odot^{\bullet}\bH^*[2].
$$
This also induces $\Z$-grading in the sheaves $\f_H$ and $\T$ on 
the supermanifold $\M$.

\sip

If $\{[e_i]\}$ is a basis in $\bH$ and $\{t^i\}$ 
are the associated linear coordinates on $\M$ as in Sect.\  2.2, then
$$
|t^i| = 2- |[e_i]|.
$$
With this choice of $\Z$-grading on $\f_{\M}$ we ensure that 
$|\Gamma|=2$ and hence $|\p|=1$, $|\delta|=1$, and $|\mu_n|=n$ for all 
the induced higher  order products on $\T$.

\sip

The {\em Euler field}\, on $M$ is, by definition, the derivation $E$ of 
$k[[t]]$ given by 
$$
E f = \frac{1}{2}|f| f, \ \ \ \ \ \forall f\in k[[t]].
$$
In coordinates,
$$
E=\frac{1}{2} \sum_i |t^i| t^i \frac{\p}{\p t^i}.  
$$
This vector field generates the scaling symmetry on $(\M, \mu_*)$ 
(cf.\ \cite{BaKo,Ma}). If we decompose
$$
\mu_n\left(\frac{\p\,}{\p t^{i_1}},\ldots,\frac{\p\,}{\p t^{i_n}}
\right) = \sum_{k} \mu_{i_1\ldots i_n}^k(t)
\frac{\p}{\p t^{k}},
$$
then, as follows from the explicit construction of $\mu_n$ given in the proof
of Theorem~4.9,
$$
E\, \mu_{i_1\ldots i_n}^k = \frac{1}{2}\left(|t^k| - |t^{i_1}| - \ldots - 
|t^{i_n}| + n\right)  \mu_{i_1\ldots i_n}^k.
$$
Note also that in the presence of identity,  $[e,E]=e$.

\bip

\paragraph{\bf 4.9.4. The perturbative proof of Theorem~A.} 
The required statement follows immediately from the graded commutative 
version of Theorem~4.9 and Sect.\ 4.9.3. \hfill $\Box$

\bip

\paragraph{\bf 4.9.5. A generalization to $G_{\infty}$-algebras.} 
In the perturbative construction of the $F_{\infty}$-functor
for dG-algebras the odd Poisson identity was used in a few 
places. For example, in the construction of $\mu_2$ the only place
where we relied on it was the cyclicity of $\Gai\cdot\Gai$,
$$
\Dg (\Gai\cdot \Gai) =0.
$$
However, a glance at the basic equation (\ref{2}) (and its higher
order analogues in Sect.\ 4.9) shows that the perturbative argument stands
if the cyclicity (and its analogues) holds only up to a homotopy.
Therefore, the generalization from dG-algebras to $G_{\infty}$-algebras
is straightforward, affecting only auxiliary tensors $A_n$.

\bip

\sip

{\small
{\em Acknowledgement.} This work was done during author's visit
to the Max Planck Institute for Mathematics in Bonn. Excellent working
conditions in the MPIM are gratefully acknowledged.
I would like 
to thank Yu.I.\ Manin for many stimulating discussions,
 and A.A.\ Voronov for  valuable communications. I am also grateful
to J.\ Stasheff for useful comments.}

\vse

{\small

{\small
\begin{tabular}{l}
Max Planck Institute for Mathematics in Bonn, and\\
Department of Mathematics, University of Glasgow\\
sm@maths.gla.ac.uk
\end{tabular}
}

\end{document}